%% file: v2-final.tex
\documentclass[noams]{moduli}

\renewcommand{\contentsname}{Contents\\{\footnotesize\normalfont(A table
of contents should normally not be included)}}

\usepackage[backend=bibtex,style=alphabetic]{biblatex}
\usepackage{mathtools}
\usepackage[all,cmtip]{xy}
\usepackage{amssymb}
\usepackage{xcolor}
\usepackage{amsmath}
\usepackage{mathrsfs}
\usepackage{enumitem}
\input{./macros.tex}

\theoremstyle{sltheoremstyle}
\newtheorem{theorem}{Theorem}[section]
\newtheorem{lemma}[theorem]{Lemma}

\newtheorem{corollary}[theorem]{Corollary}

\newtheorem{proposition}[theorem]{Proposition}

\newtheorem*{corollary-non}{Corollary}
\newtheorem*{lemma-non}{Lemma}
\newtheorem*{theorem-non}{Theorem}
\newtheorem*{proposition-non}{Proposition}
\newtheorem*{condition-non}{Condition}
\newtheorem*{conditions-non}{Conditions}

\theoremstyle{definition}

\newtheorem{remark}[theorem]{Remark}
\newtheorem{remarks}[theorem]{Remarks}
\newtheorem{definition}[theorem]{Definition}

\begin{document}

\title{Hyperbolic geometry and real moduli of five points on the line}
\author{Olivier de Gaay Fortman}
\email{
a.o.d.degaayfortman@uu.nl
}
\address{
 Department of Mathematics,
Utrecht University,
       Budapestlaan 6, 3584 CD Utrecht, The Netherlands
}
%
%
\keywords{Binary quintics, moduli spaces, period maps, real algebraic geometry, hyperbolic geometry. \\ \emph{Date:} \today.}

\begin{abstract}

We show that each connected component of the moduli space of smooth real binary quintics is isomorphic to an open subset of an arithmetic quotient of the real hyperbolic plane. Moreover, our main result says that the induced metric on this moduli space extends to a complete real hyperbolic orbifold structure on the space of stable real binary quintics. This turns the moduli space of stable real binary quintics into the quotient of the real hyperbolic plane by an explicit non-arithmetic triangle group. 
\vspace*{-6.5pt}
\end{abstract}

\maketitle


\section{Introduction} \label{realbinaryintroduction}

For interesting classes of complex varieties, there is a period map that identifies the moduli space with an open subset of an arithmetic quotient of a hermitian symmetric domain. Classical examples include abelian varieties, K3-surfaces, and configurations of points on the line. To study moduli of real algebraic varieties, several authors analyzed the equivariance of the complex period map with respect to the action of complex conjugation on cohomology \cite{kharlamov-surfaces, nikulin-intsymforms, MR739089, seppalasilhol, MR1668340, MR1623240, aperyyoshida-pentagonalstructure, MR1795406}. 
An important difference between the complex and real case is that moduli spaces of smooth real varieties are often not connected. This implies that a real period map has to be defined on each connected component of the moduli space separately; in favourable cases, this defines an isomorphism between any such a component and the quotient of a Riemannian manifold by a discrete group of isometries, see e.g.~\cite{grossharris, realACTsurfaces, MR2848997, heckman2016hyperbolic}. 

To salvage the non-connectedness of the real moduli space, one can sometimes define a slightly larger moduli space by allowing mild singularities. The idea is that in such a larger space, the smooth varieties of one topological type now do deform into smooth varieties of another topological type, making the moduli space \emph{connected}. 
In their beautiful paper \cite{realACTsurfaces}, Allcock, Carlson and Toledo showed that for cubic surfaces, the real period maps defined on the various connected components of the moduli space of smooth surfaces extend to a global period map, defined on the moduli space of stable real cubic surfaces. In this way, they identified the latter with a single non-arithmetic real hyperbolic quotient. 
They proved analogous results for moduli of stable binary sextics, and stable binary sextics with a double root at infinity  \cite{realACTnonarithmetic,realACTbinarysextics}. 

It turns out that binary quintics provide a new example of this phenomenon. 
Let $X \cong \bb A^6_\RR$ be the real algebraic variety that parametrizes homogeneous polynomials $F \in \RR[x,y]$ of degree five. Let 
$X_0 \subset X$ parametrize polynomials with distinct roots, and $X_s \subset X$ polynomials with roots of multiplicity at most two (i.e.~stable in the sense of geometric invariant theory). The principal goal of this paper is to study the geometry of the moduli space of stable \textit{real binary quintics} 
$$
 \ca M_s(\RR) \coloneqq \GL_2(\RR) \setminus X_s(\RR)  \supset   \GL_2(\RR) \setminus X_0(\RR) \eqqcolon \ca M_0(\RR).  
$$
Let $P_s \subset \PP^1(\CC)^5$ be the set of five-tuples $(x_1, \dotsc, x_5)$ such that no three $x_i \in \PP^1(\CC)$ coincide (cf.~\cite{MR0437531}), and let $P_0 \subset P_s$ be the subset of five-tuples whose coordinates are distinct. These spaces are naturally acted upon by $\mf S_5$, the symmetric group on five letters. Moreover, complex conjugation $\sigma \colon \PP^1(\CC)^5 \to \PP^1(\CC)^5$ 
induces an 
anti-holomorphic involution $\sigma \colon P_s/\mf S_5 \to P_s/\mf S_5$ that preserves $P_0/\mf S_5$. Let $(P_0/\mf S_5)(\RR)$ and  $(P_s/\mf S_5)(\RR)$ denote the respective fixed loci. Then
$$
\ca M_0(\RR) \cong \PGL_2(\RR) \setminus (P_0/\mf S_5)(\RR) \white \white \text{ and } \white \white \ca M_s(\RR) \cong \PGL_2(\RR) \setminus (P_s/\mf S_5)(\RR). 
$$ 
In other words, $\ca M_0(\RR)$ is the space of subsets $S \subset \PP^1(\CC)$ of cardinality $|S| = 5$ stable under complex conjugation modulo real projective transformations; in $\ca M_s(\RR)$ one or two pairs of points are allowed to collapse. 

By Deligne--Mostow theory, the coarse moduli space $\ca M_s(\CC)= \GL_2(\CC) \sm X_s(\CC)$ of stable complex binary quintics has a complex hyperbolic orbifold structure. 
Indeed, for five distinct points $u_1, \dotsc, u_5 \in \bb A^1(\CC) \subset \PP^1(\CC)$, the projective model of the normalization of the affine curve $z^5 = (x-u_1)^2 \cdots (x-u_5)^2$ is a smooth curve $C$ of genus six; this curve $C$ carries an automorphism of order five that induces an automorphism on the space $\rm H^0(C, \Omega^1_C)$ of holomorphic one-forms on $C$ whose $e^{2 \pi i/5}$-eigenspace defines a line in the corresponding eigenspace in $\rm H^1(C(\CC),\CC)$. This line is negative for a natural hermitian form, and hence one can associate to $\set{u_1, \dotsc, u_5}$ a point in a certain two-dimensional complex ball quotient $P\Gamma \sm \CCH^2$. This construction was already known to Shimura, see \cite{Shimura1963ONAF, shimuratranscendental}. 
By varying the subset $\set{u_1, \dotsc, u_5}$ of points on $\PP^1(\CC)$, or rather the associated complex binary quintic, one obtains a period map $\ca M_0(\CC) = \GL_2(\CC) \sm X_0(\CC) \to P\Gamma \sm \CCH^2$ 
(see Section \ref{section:periodmap} for details). Results of Deligne and Mostow \cite{DeligneMostow} imply that this period map 
extends to 
an isomorphism of complex analytic spaces $
\ca M_s(\CC) = \GL_2(\CC) \setminus X_s(\CC) \xrightarrow{\sim} P\Gamma \setminus \CC H^2$, see Theorem \ref{th:delignemostow}. Since strictly stable quintics correspond to points in a hyperplane arrangement $\mr H \subset \CC H^2$, see Proposition \ref{prop:stableperiodshyperplane}, one thus obtains an isomorphism
\begin{align} \label{equation:periodmap-intro}
\ca M_0(\CC) = \GL_2(\CC) \setminus X_0(\CC) \xrightarrow{\sim} P\Gamma \setminus \left(\CC H^2 - \mr H \right).
\end{align}
By investigating the equivariance of the period map with respect to suitable anti-holomorphic involutions $\alpha_j \colon \CC H^2 \to \CC H^2,$ we obtain the following real analogue.


\begin{theorem} \label{th:theorem01}
For $j \in \{0,1,2\}$, let $\mr M_{j}$ be the connected component of $\ca M_0(\RR)$ parametrizing $\Gal(\CC/\RR)$-stable subsets $S \subset \PP^1(\CC)$ with $2j$ complex and $5 - 2j$ real points. The period map induces an isomorphism of real analytic orbifolds 
\begin{align} \label{iso:smoothcase}
\mr M_j \xrightarrow{\sim} P \Gamma_j \setminus \left(\RR H^2 - \mr H_j \right).
\end{align}
Here $\RR H^2$ is the real hyperbolic plane, $\mr H_j$ a union of geodesic subspaces in $ \RR H^2$, and $P\Gamma_j$ an arithmetic lattice in $\rm{PO}(2,1)$. Moreover, the lattices $P\Gamma_j$ are projective orthogonal groups attached to explicit quadratic forms over $\ZZ[\lambda]$,  $\lambda = \zeta_5 + \zeta_5^{-1} = (\sqrt 5 - 1)/2$, see equation \eqref{eq:explicitquadraticforms}. 
\end{theorem}
In particular, 
Theorem \ref{th:theorem01} endows each connected component $\mr M_j \subset \ca M_0(\RR)$ with a hyperbolic metric. Since one can deform the topological type of a $\Gal(\CC/\RR)$-stable five-element subset of $\PP^1(\CC)$ by allowing two points to collide, the compactification $\ca M_s(\RR) \supset \ca M_0(\RR)$ is connected. One may wonder whether the metrics on the components $\mr M_j$ extend to a metric on the whole of $\ca M_s(\RR)$. If so, what does the resulting space look like at the boundary? Our main result answers these questions in the following way.

\begin{theorem} \label{th:theorem02}
There exists a complete hyperbolic metric on $\ca M_s(\RR)$ that restricts to the metrics on $\mr M_j$ induced by (\ref{iso:smoothcase}). 
Let $\overline{\mr{M}}_\RR$ denote the resulting metric space, and define $\Gamma_{3,5,10}$ as the group
\begin{align}\label{PGAMMAR}
\Gamma_{3,5,10} =  \langle a, b, c \mid a^2 = b^2 = c^2 = (ab)^3 = (ac)^5 = (bc)^{10} = 1 \rangle.
\end{align} There exists an open embedding $\coprod_j P \Gamma_j \setminus \left(\RR H^2 - \mr H_j \right) \hookrightarrow \Gamma_{3,5,10}\setminus \RR H^2$ and an isometry 
\begin{align} \label{isometry}
\overline{\mr{M}}_\RR \cong \Gamma_{3,5,10} \setminus \RR H^2
\end{align} 
that extend the real analytic orbifold isomorphisms (\ref{iso:smoothcase}) in Theorem \ref{th:theorem01}. 
In particular, $\overline{\mr{M}}_\RR$ is isometric to the hyperbolic triangle $\Delta_{3,5,10}$ of angles $\pi/3, \pi /5, \pi/10$, see Figure~\ref{fig:triangle} below.
\end{theorem}


\hspace*{-6cm}\includegraphics[scale=0.24]{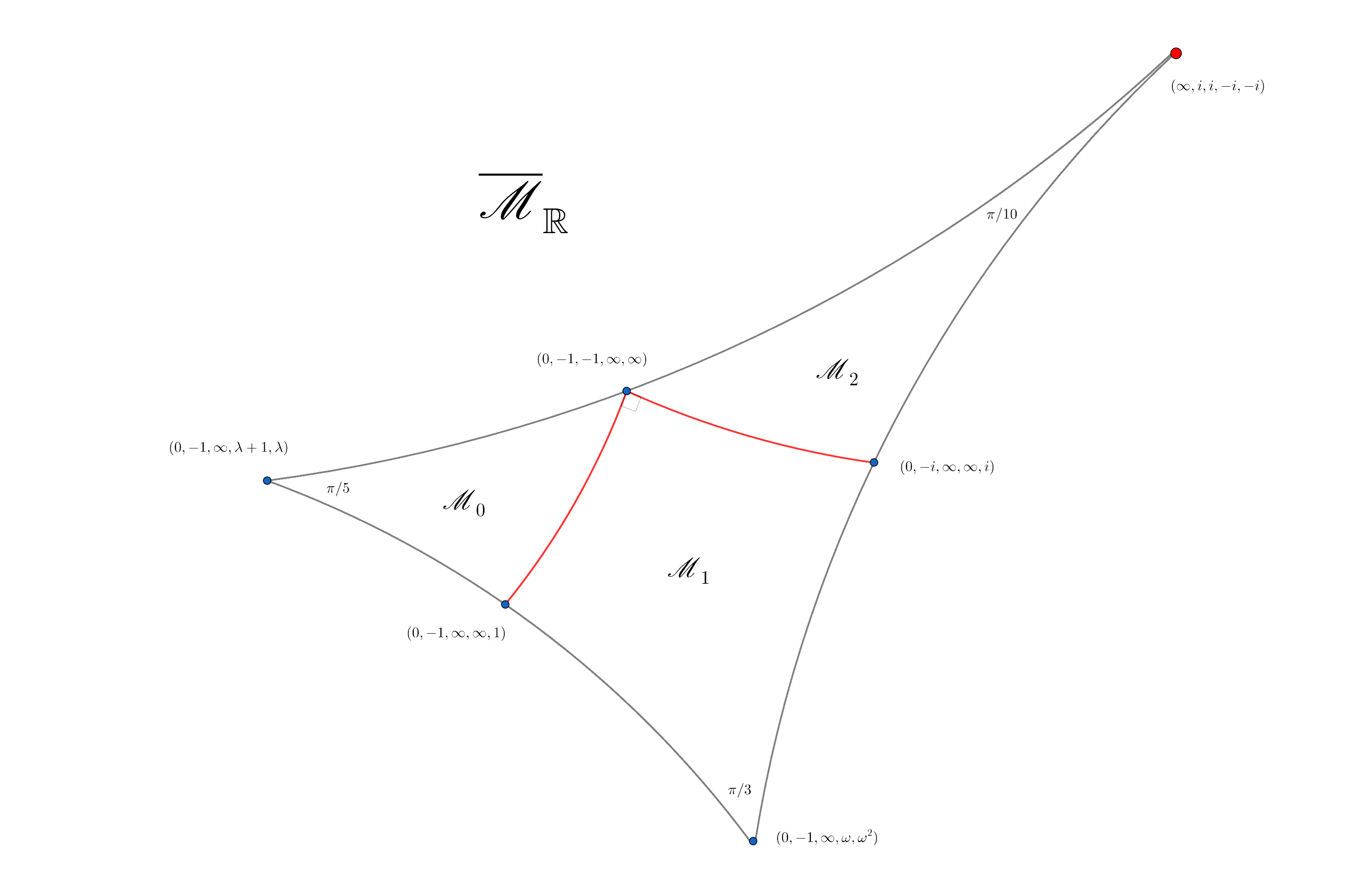}     \label{fig:triangle}
\noindent
\emph{Figure 1: 
The moduli space of stable real binary quintics as the hyperbolic triangle $\Delta_{3,5,10} \subset \RR H^2$. Here $\lambda = \zeta_5 + \zeta_5^{-1} = (\sqrt 5 - 1)/2$ and $\omega = \zeta_3$ (where $\zeta_n = e^{2 \pi i/n} \in \CC$ for $n \in \ZZ_{\geq 3}$).}
\\
\\
Note that the closure $$\overline{\mr M}_0 \subset \overline{\mr M}_\RR$$ of $\mr M_0$ in $\overline{\mr M}_\RR$ is the moduli space of stable configurations of five real points on $\PP^1_\RR$. This moduli space was investigated by Apéry and Yoshida in \cite{aperyyoshida-pentagonalstructure}, who proved that $\overline{\mr M}_0$ is the hyperbolic triangle of angles $\pi/2, \pi/4$ and $\pi/5$. From this, together with the knowledge of the angles of $\overline{\mr M}_\RR$ and the fact that the two hyperplanes in Figure \ref{fig:triangle} intersect orthogonally, one can deduce the remaining angles of the closures $\overline{\mr M}_j \subset \overline{\mr M}_\RR$ of the subsets $\mr M_j \subset \overline{\mr M}_\RR$ ($j \in \set{0,1,2}$). 
\begin{theorem} \label{theorem:angles}
Consider Figure \ref{fig:triangle}. For $j = 0,1,2$, let $\overline{\mr M}_j \subset \overline{\mr M}_\RR$ be the closure of $\mr M_j \subset \overline{\mr M}_\RR$.
\begin{enumerate}
\item
The angle of $\overline{\mr M}_0$ at $(0,-1,\infty,\infty,1)$ is $\pi/2$, and its angle at $(0,-1,-1,\infty,\infty)$ is $\pi/4$. 
\item The angle of $\overline{\mr M}_1$ at $(0,-1, \infty,\infty,1)$ is $\pi/2$, and its angle at $(0,-i,\infty,\infty,i)$ is $\pi/2$. \item The angle of $\overline{\mr M}_2$ at $(0,-1,-1,\infty,\infty)$ is $\pi/4$, and its angle at  $(0,-i,\infty,\infty,i)$ is $\pi/2$. 
\end{enumerate}
\end{theorem}

\begin{remarks} \label{rem:introremarks}
\begin{enumerate}
\item 
 \label{remark:takeuchi}
The lattice $\Gamma_{3,5,10} \subset \textnormal{PO}(2,1)$ is \textit{non-arithmetic}, see \cite{takeuchi}. 
\item
The topological space $\ca M_s(\RR)$ underlies two topological orbifold structures: the natural orbifold structure on $ \GL_2(\RR) \setminus X_s(\RR)$ and the orbifold structure on $\overline{\mr M}_\RR$ induced by (\ref{isometry}). 
These orbifold structures only differ at one point of $\ca M_s(\RR)$, which is $(\infty, i, i, -i, -i)$ (see Figure~\ref{fig:triangle}). The stabilizer group of $\ca M_s(\RR)$ at $(\infty, i, i, -i, -i)$ is isomorphic to $\ZZ/2$, whereas the stabilizer group of $\overline{\mr M}_\RR$ at $(\infty, i, i, -i, -i)$ is isomorphic to the dihedral group of order $20$.
\item Important ingredients in the proof of Theorem \ref{th:theorem02} 
are: 
\begin{enumerate}
\item \label{item:a}
the fact that under the complex period map \eqref{equation:periodmap-intro}, moduli of singular binary quintics correspond to points in the quotient of a certain hyperplane arrangement $\mr H \subset \CC H^2$; and 
\item \label{item:b} the fact that the hyperplane arrangement $\mr H \subset \CC H^2$ is an \textit{orthogonal arrangement} in the sense of \cite{orthogonalarrangements}. 
\end{enumerate}
We prove (a) in Proposition \ref{prop:stableperiodshyperplane}, and (b) holds by \cite[Theorem 2.5 \& Example 2.12]{degaayfortman-nonarithmetic}. 
\end{enumerate}
\end{remarks}
\begin{remark}
Let $PX_0(\CC)$ denote the space of $\CC^\ast$-equivalence classes of smooth complex binary quintics $F \in \CC[x,y]$. The natural map $P_0 \to PX_0(\CC)$ 
induces a $\PGL_2(\CC)$-equivariant isomorphism $P_0/\mf S_5 \xrightarrow{\sim} PX_0(\CC)$, and the quotient $\ca M_0(\CC) = \PGL_2(\CC) \sm PX_0(\CC)$ 
is the moduli space of smooth complex binary quintics. 
It turns out that neither $\pi_1 \left( PX_0(\CC)  \right)$ nor $\pi_1^\textnormal{orb}\left( \ca M_0(\CC) \right)$ is a lattice in any Lie group with finitely many connected components. In view of \cite[Theorem 1.2]{orthogonalarrangements}, this follows from the isomorphism $\ca M_0(\CC) \cong P\Gamma \setminus \left(\CC H^2 - \mr H \right)$ (see equation \eqref{equation:periodmap-intro} above) and the orthogonality of the hyperplane arrangement $\mr H \subset \CC H^2$ (see Remark \ref{rem:introremarks}.(iii).(b) above). 
\end{remark}


\subsection{Overview of this paper} In Section \ref{section:monodromy}, we recall known results on families of quintic covers of the complex projective line, branched along a binary quintic. We consider moduli of complex binary quintics in Section \ref{complexball}. In particular, we show that Deligne--Mostow theory provides an isomorphism between the space of stable complex binary quintics and an arithmetic ball quotient. 
In Section \ref{modrealbinquin}, we prove that 
moduli of stable real binary quintics are in one-to-one correspondence with points in the real hyperbolic quotient space $P\Gamma_\RR\setminus \RR H^2$ defined by a lattice $P\Gamma_\RR \subset \rm{PO}(2,1)$. We calculate $P\Gamma_\RR$ in Section \ref{space:hyperbolictriangle}: it is conjugate to the lattice $\Gamma_{3,5,10}$ defined in \eqref{PGAMMAR}. In Section \ref{sec:monodromy}, we study monodromy groups of moduli spaces of smooth binary quintics over $\CC$ and over $\RR$, and prove Theorems \ref{th:theorem02} and \ref{theorem:angles}. 
 In Section \ref{sec:announced}, 
we use \cite[Theorem 1.8]{degaayfortman-nonarithmetic} and the main results of this paper to 
provide an explicit sequence $\set{\Gamma_n \subset \rm{PO}(n,1)}_{n \geq 2}$ of non-arithmetic lattices $\Gamma_n$, with $\Gamma_2 = P\Gamma_\RR$. 

\section{Notation} \label{section:binarynotation}

Let $K$ be the cyclotomic field $ \QQ(\zeta)$, where $$\zeta \coloneqq \zeta_5 = e^{2 \pi i /5} \in \CC.$$ The ring of integers $\OO_K$ of $K$ is $\ZZ[\zeta]$, see e.g.~\cite[Chapter I, Proposition 10.2]{Neukirch}. Let $\mu_K \subset \OO_K^\ast$ be the torsion subgroup of the unit group $\OO_K^\ast$, and recall that $\mu_K$ is cyclic of order ten, generated by $-\zeta$. Define an involution $\rho \colon K \to K$ by $\rho(\zeta) = \zeta^{-1}$, and let $F = K^\rho$ be the maximal totally real subfield of $K$. Recall that $F$ is generated over $\QQ$ by  the element
$$\lambda \coloneqq \zeta + \zeta^{-1} = (\sqrt 5 - 1)/2.$$
Define 
\[
\eta = \zeta^2 - \zeta^{-2} \in \OO_K,
\] 
and consider the different ideal $\mf D_K \subset \OO_K$, see e.g.~\cite[Chapter III, Section 2, Definition 2.1]{Neukirch} or \cite[Chapter III, Section 3]{localfields}. 
We have 
$
(\zeta - \zeta^{-1}) \cdot (\zeta^2 - \zeta^{-2})  \cdot (\zeta^3 - \zeta^{-3}) \cdot (\zeta^4 - \zeta^{-4}) = 5,
$
hence $
(\eta)^4 = (5)
$
as ideals of $\OO_K$. This implies (see e.g.~\cite[Chapter III, Theorem 2.6]{Neukirch}) that
\[
\mf D_K = \left(5/ \eta \right) = (\eta)^3. 
\]

\section{Preliminaries on quintic covers of the projective line branched along five points}\label{section:monodromy}

In this section, we recollect some known results 
on quintic covers $C \to \PP^1_\CC$ ramified along five points with local monodromy $\exp(4 \pi i / 5)$ around each point. 
Some of these results are well-known but hard to find in the literature; we state and prove those for convenience of the reader. Other results of this section are available in the literature, but we formulate them in a different manner. 


Recall from the introduction that $X \cong \bb A^6_\RR$ is the real algebraic variety that parametrizes homogeneous polynomials of degree five, $X_0$ the subvariety of polynomials with distinct roots, and $X_s \subset X$ the subvariety of polynomials with roots of multiplicity at most two, i.e.~non-zero polynomials whose class in the associated projective space is stable in the sense of geometric invariant theory \cite{GIT} for the action of $\SL_{2, \RR}$ on it. 





\subsection{Refined Hodge numbers of a quintic cover of $\PP^1$ branched along five points} \label{jacofcyc} 

Let $F \in X_0(\CC)$ be a smooth binary quintic; thus $F = F(x,y) \in \CC[x,y]$ is a homogeneous polynomial of degree five whose zeros in $\PP^1(\CC)$ all have multiplicity one. Let $t_1, \dotsc, t_5 \in \PP^1(\CC)$ be the zeros of $F$, with $t_j = [u_j \colon v_j]$ in homogeneous coordinates of $\PP^1(\CC)$, for $j = 1, \dotsc, 5$. Let
\begin{align} \label{equation:associatedcycliccover}
C_F \to \PP^1_\CC
\end{align}
be the cyclic quintic cover with branch points $t_1, \dotsc, t_5 \in \PP^1(\CC)$ and local monodromy 
\begin{align} \label{equation:localmonodromy}
\exp(2 \pi i \cdot 2/5) \in \mu_5
\end{align}
around $t_j$ for each $j \in \set{1, \dotsc, 5}$. If the points $t_j = [u_j \colon 1]$ are all in $\bb A^1(\CC)$, the cover \eqref{equation:associatedcycliccover} is in affine coordinates given by the normalization of the curve defined by the equation
\begin{align*}
z^5 = (x - u_1)^2 \cdot (x - u_2)^2 \cdots (x-u_5)^2,
\end{align*}
with $\zeta \in \mu_5 \subset \CC$ acting by $(x,z) \mapsto (x, \zeta \cdot z)$. We have $g = 6$ for the genus $g = g(C_F)$ of the curve $C_F$. Let $JC_F$ be the Jacobian of the curve $C_F$, so that
\[
JC_F(\CC)= \rm H^1(C_F(\CC), \OO_{C_F})/\rm H^1(C_F(\CC), \ZZ)(1),
\] with weight $-1$ Hodge decomposition 
\begin{align}\label{eq:hodgedecomp}
\rm H_1(JC_F(\CC), \CC) = \rm H^1(C_F(\CC), \ZZ)(1) \otimes \CC = \rm H^{-1,0}(JC_F) \oplus \rm H^{0,-1}(JC_F).
\end{align}
Note that $\rm H^{0,-1}(JC_F)$ is naturally isomorphic to the space $\rm H^0(C_F, \Omega^1)$ of global holomorphic differentials on the curve $C_F$. The order five automorphism $$\zeta \colon C_F \xrightarrow{\sim}C_F$$ defined above induces an embedding of rings
\[
\varphi \colon \ZZ[\zeta] \to \End(JC_F), \quad \varphi(\zeta) = (\zeta^{-1})^\ast,
\]
which is compatible with the Hodge decomposition \eqref{eq:hodgedecomp}. For $k \in \set{1,2,3,4}$, define 
\[
\rm H^{0,-1}(JC_F)_{\zeta^k} = \set{x \in \rm H^{0,-1}(JC_F) \mid \varphi(\zeta) = \zeta^k} \subset \rm H^{0,-1}(JC_F) = \rm H^{1,0}(C_F) = \rm H^0(C_F, \Omega^1),
\]
and define $\rm H^{-1,0}(JC_F)_{\zeta^k} \subset \rm H^{-1,0}(JC_F)$ in a similar way. 


It is well-known how to calculate the refined Hodge numbers $h^{0,-1}(JC_F)_{\zeta^k} = \dim \rm H^{0,-1}(JC_F)_{\zeta^k}$ and $h^{-1,0}(JC_F)_{\zeta^k}= \dim \rm H^{-1,0}(JC_F)_{\zeta^k}$ for $k = 1,2,3,4$; the result is as follows. 

\begin{lemma} \label{lemma:refinedhodge}
Let $F \in X_0(\CC)$ be a smooth binary quintic, and let $JC_F$ be the Jacobian of the cyclic cover $C_F \to \PP^1(\CC)$ associated to $F$ as in \eqref{equation:associatedcycliccover}. 
One has the following refined Hodge numbers: 
\begin{align*}
    h^{0,-1}(JC_F)_{\zeta} &= 1, \white h^{0,-1}(JC_F)_{\zeta^2} = 3, \white 
    h^{0,-1}(JC_F)_{\zeta^3} = 0, \white 
    h^{0,-1}(JC_F)_{\zeta^{4}} = 2, \\
        h^{-1,0}(JC_F)_{\zeta} &= 2, \white h^{-1,0}(JC_F)_{\zeta^2} = 0, \white 
    h^{-1,0}(JC_F)_{\zeta^3} = 3, \white 
    h^{-1,0}(JC_F)_{\zeta^{4}} = 1.
\end{align*}
\end{lemma}
\begin{proof}
This follows from the Hurwitz-Chevalley--Weil formula, see \cite[Proposition 5.9]{Moonen2011TheTL}. Alternatively, see \cite[Lemma 4.2]{looijenga-lauricella}.
\end{proof}

\subsection{The hermitian lattice} \label{section:hermitianlattice}

We fix a smooth binary quintic $F_0 \in X_0(\CC)$. Let 
$
C \to \PP^1_\CC
$ 
be the cyclic cover of $\PP^1_\CC$ associated to $F_0$ as in \eqref{equation:associatedcycliccover}, and consider the Jacobian variety with $\OO_K$-action
\[
\left(JC ,\quad \varphi \colon \OO_K = \ZZ[\zeta] \to \End(A)\right).
\]
See Section \ref{jacofcyc} above. Define $\Lambda$ as the free $\OO_K$-module
\[
\Lambda =\rm H^1(C(\CC), \ZZ). 
\]
The canonical principal polarization of the Jacobian $JC$ is given by a symplectic pairing
\[
E \colon \Lambda \times \Lambda = \rm H^1(C(\CC), \ZZ) \times \rm H^1(C(\CC), \ZZ) \to \rm H^2(C(\CC),\ZZ) = \ZZ 
\]
which satisfies $E(\varphi(a)x,y) = E(x, \varphi(\rho(a))y)$ for each $a \in \OO_K$ and $x,y \in \Lambda$, where $\rho$ is the automorphism $K \xrightarrow{\sim} K, \zeta \mapsto \zeta^{-1}$, see Section \ref{section:binarynotation}. 
Consider the different ideal $\mf D_K \subset \OO_K$, and define a skew-hermitian form $T$ on $\Lambda$ as follows:
\[
T \colon \Lambda \times \Lambda \to \mf D_K^{-1},  \quad T(x,y) = \frac{1}{5}\sum_{j = 0}^{4}\zeta^jE\left( x, \varphi(\zeta)^j y \right).
\]
By \cite[Example 11.2.2]{degaayfortman-nonarithmetic}, this is the skew-hermitian form corresponding to $E$ via \cite[Lemma 11.1]{degaayfortman-nonarithmetic}. Let $\eta \in \OO_K$ be the purely imaginary element
$$
\eta = \zeta^2 - \zeta^{-2} \in \OO_K,
$$
and remark that the different ideal $\mf D_K \subset \OO_K$ is generated by $5/\eta$, see Section \ref{section:binarynotation}. We obtain a hermitian form on the free $\OO_K$-module $\Lambda$ as follows:
\begin{equation} \label{eq:hermitianformonbinaryquinticlattice}
    \mf h \colon \Lambda \times \Lambda \to \OO_K, \white
\mf h(x,y) = \frac{5}{\eta} \cdot T(x,y) = \frac{1}{\zeta^2 - \zeta^{-2}} \cdot \sum_{j = 0}^{4}\zeta^jE\left( x, \varphi(\zeta)^j y \right). 
\end{equation}
By \cite[Lemma 11.1]{degaayfortman-nonarithmetic}, we have that $(\Lambda, \mf h)$ is unimodular as a hermitian lattice over $\OO_K$, because $(\Lambda, E)$ is unimodular as an alternating lattice over $\ZZ$. 

\subsection{The signature} \label{subsection:signature}

Define a CM type $\Phi \subset \Hom(K,\CC)$ as
\begin{align} \label{equation:CMtype}
\Phi = \set{\sigma_1, \sigma_2 \colon K \to \CC}, \quad \sigma_k(\zeta) = \zeta^k \quad \text{for}\quad k = 1,2. 
\end{align}
Observe that, for $k = 1,2$, we have
$$
\Lambda \otimes_{\OO_K, \sigma_k}\CC= \left(\Lambda \otimes_\ZZ \CC\right)_{\zeta^k} = \rm H^1(C(\CC),\CC)_{\zeta^k}.
$$

\begin{lemma} \label{lemma:signature}
For $k \in \set{1,2}$, consider the hermitian form 
$$
\mf{h}^{\sigma_k}  \colon \rm H^1(C(\CC),\CC)_{\zeta^k} \times \rm H^1(C(\CC),\CC)_{\zeta^k} \to \CC.
$$
Then $\tn{sign}(\mf{h}^{\sigma_1})  = (2,1)$ and $\tn{sign}(\mf h^{\sigma_2}) = (3,0)$, where $\tn{sign}(\mf h^{\sigma_k})$ is the signature of $\mf{h}^{\sigma_k}$. 
\end{lemma}

\begin{proof}
Write $\Lambda_\CC = \Lambda \otimes_\ZZ \CC = \rm H^1(C(\CC),\CC)$. 
For each embedding $\phi \colon K \to \CC$, the restriction of the hermitian form $$\phi\left(5/\eta \right) \cdot E_\CC(x, \bar y) \colon \Lambda_\CC \times \Lambda_\CC \to \CC$$
to $(\Lambda_\CC)_{\phi} \subset \Lambda_\CC$ coincides with 
$
\mf{h}^{\phi}$  
by \cite[Lemma 11.3]{degaayfortman-nonarithmetic}. Moreover, the hermitian form
\[
i \cdot E_\CC(x, \bar y) \colon \rm H^1(C(\CC),\CC) \times \rm H^1(C(\CC),\CC) \to \CC
\] is positive definite on 
$
\rm H^0(C, \Omega^1) = \rm H^{1,0}(C) = \rm H^{0,-1}(JC)
$
and negative definite on $\rm H^{-1,0}(JC) = \rm H^{0,1}(C)$, see \cite[Théorème 6.32]{voisin}. 
As $\Im(\sigma_1(\eta)) > 0$ and $\Im(\sigma_2(\eta)) < 0$, we have $\Im(\sigma_1(5/\eta)) < 0$ and $\Im(\sigma_2(5/\eta)) > 0$. Consequently, the hermitian form 
$$
\sigma_1(5/\eta)\cdot E_\CC(x, \bar y) = \mf{h}^{\sigma_1}(x,y) \colon \rm H^1(C(\CC),\CC)_\zeta \times \rm H^1(C(\CC),\CC)_\zeta  \to \CC
$$
is negative definite on $\rm H^{0,-1}(JC)_\zeta$ and positive definite on $\rm H^{-1,0}(JC)_{\zeta}$, so that $
\tn{sign}(\mf{h}^{\sigma_1}) = (2,1)$ by Lemma \ref{lemma:refinedhodge}. Similarly, the hermitian form 
$$
\sigma_2(5/\eta)\cdot E_\CC(x, \bar y) = \mf{h}^{\sigma_2}(x,y) \colon \rm H^1(C(\CC),\CC)_{\zeta^2} \times \rm H^1(C(\CC),\CC)_{\zeta^2}  \to \CC
$$
is positive definite on $\rm H^{0,-1}(JC)_{\zeta^2}$ and negative definite on $\rm H^{-1,0}(JC)_{\zeta^2}$. Hence, using Lemma \ref{lemma:refinedhodge} again, we conclude that $\tn{sign}(\mf h^{\sigma_2}) = (3,0)$.
\end{proof}


\subsection{The monodromy representation} \label{sec:monodromy:two-three} Consider the real algebraic variety $X_0$ introduced in Section \ref{realbinaryintroduction}. 
Let $D \subset \GL_2(\CC)$ be the subgroup $D = \set{\zeta^j\cdot \Id } \subset \GL_2(\CC)$ of scalar matrices of the form $\zeta^j \cdot \Id$ with $j \in \set{0,1,2,3,4}$ and $\Id \in \GL_2(\CC)$ the identity two by two matrix. Define
\begin{align} \label{def:GC}
G(\CC) = \GL_2(\CC)/D.  
\end{align}
The group $G(\CC)$ acts from the left on $X_0(\CC)$ in the following way: if $$F = F(x,y) \in \CC[x,y]$$ is a binary quintic, we may view $F$ as a function $\CC^2 \to \CC$, and define $g \cdot F = F(g^{-1})$ for $g \in G(\CC)$. This gives a canonical isomorphism of complex analytic orbifolds
$$
\ca M_0(\CC) = G(\CC) \setminus X_0(\CC),
$$
where $\ca M_0$ is the moduli stack of smooth binary quintics. 
Define 
\begin{align} \label{universalcycliccover}
g \colon \mr C \to X_0 
\end{align} 
as the universal family of cyclic quintic covers $C \to \PP^1$ ramified along a smooth binary quintic $\{F = 0\} \subset \PP^1$ with local monodromy given by \eqref{equation:localmonodromy}, and let 
 $$
 J = J(\mr C) \to X_0$$ 
 be the relative Jacobian of $\mr C/X_0$. By Lemma \ref{lemma:signature}, $J$ is a polarized abelian scheme of relative dimension six over $X_0$, equipped with $\OO_K$-action of signature $\{(2,1), (3,0)\}$ with respect to $\Phi = \{\sigma_1, \sigma_2\}$. Let $\bb V = R^1g_\ast\ZZ$ be the local system of hermitian $\OO_K$-modules underlying the abelian scheme $J/X_0$. Attached to $\bb V$ and the base point $F_0 \in X_0(\CC)$, we have a representation 
\begin{align} \label{eq:monodromy-nonprojective}
\rho_\Gamma \colon \pi_1(X_0(\CC), F_0) \to \Gamma \coloneqq \Aut_{\OO_K}(\Lambda, \mf h)
\end{align} whose composition with the quotient map $\Gamma \to P\Gamma = \Gamma / \mu_K$ defines a homomorphism 
\begin{equation} \label{eq:monodromy}
P(\rho_\Gamma) \colon \pi_1(X_0(\CC), F_0) \to P\Gamma.
\end{equation}

\subsection{Marked binary quintics and points on the projective line} \label{sec:marked} Let $F  \in X_0(\CC)$ and consider the hypersurface
$
 Z_F \coloneqq \{F = 0\} \subset \PP^1_\CC$. 
  A \textit{marking} of $F$ is an ordering $m \colon Z_F(\CC) \xrightarrow{\sim} \set{1, 2, 3,4,5}$ of the five-element set $Z_F(\CC) = \set{x \in \PP^1(\CC) \mid F(x) = 0}$. \begin{remark} Let $F \in X_0(\CC)$. To give a marking of $F$ is to give an isomorphism of rings $\rm H^0(Z_F(\CC), \ZZ) \xrightarrow{\sim} \ZZ^5$. 
 Moreover, sending a ring automorphism of $\ZZ^5$ to the induced permutation of the canonical basis of $\set{e_1, \dotsc, e_5} \subset \ZZ^5$ defines a group isomorphism $\Aut_{\rm{ring}}(\ZZ^5) \cong \mf S_5$.\end{remark} 
Define $\ca N_0$ as the set of marked smooth complex binary quintics $(F, m)$, and consider the map 
that forgets the marking:
\begin{align}\label{markedquintics}
\ca N_0 \to X_0(\CC), \quad \quad (F,m) \mapsto F.
\end{align}
To provide $\ca N_0$ with a complex manifold structure such that \eqref{markedquintics} is a finite covering map, we define $\ca N_0$ in a slightly different but equivalent way. Define algebraic subgroup
\[
T\coloneqq \set{ (\lambda_1, \dotsc, \lambda_5) \in (\CC^\ast)^5 \mid \lambda_1 \cdots \lambda_5 = 1} \subset (\CC^\ast)^5. 
\]
The group $T$ acts freely on $(\CC^2 - \set{0})^5$ in the usual way, i.e.~by letting $(\lambda_1, \dotsc, \lambda_5) \in T$ act as 
\[
\left( (a_1, b_1), \dotsc, (a_5, b_5) \right) \mapsto 
\left( \lambda_1 \cdot (a_1, b_1), \dotsc, \lambda_5 \cdot (a_5, b_5) \right). 
\]
As $\dim(T) = 4$, the quotient variety
$
\ca N \coloneqq (\CC^2 - \set{0})^5 / T
$
is a complex manifold of dimension six, equipped with a holomorphic principal $\CC^\ast$-bundle map 
\begin{align}\label{align:f}
f \colon \ca N \to \PP^1(\CC)^5, \quad \quad \left(v_1, \dotsc, v_5 \right) \mapsto \left([v_1], \dotsc, [v_5]\right),
\end{align}
where $\CC^\ast$ acts on $\ca N$ by $\lambda \cdot (v_1, v_2, \dotsc, v_5) = (\lambda v_1, v_2, \dotsc, v_5) =  (v_1, \lambda v_2, \dotsc, v_5) = \cdots = (v_1, v_2, \dotsc, \lambda v_5) \in \ca N$ for $\lambda \in \CC^\ast$ and $(v_1, v_2, \dotsc, v_5) \in \ca N$. For $i,j \in \set{1, \dotsc, 5}$ with $i<j$, define $\overline{\Delta}_{ij} = \set{(x_1, \dotsc, x_5) \in \PP^1(\CC)^5 \mid x_i = x_j}$, and let 
$$\Delta_{ij} \coloneqq f^{-1}(\overline{\Delta}_{ij}),  
$$
where $f$ is the principal bundle map \eqref{align:f}. Since $\overline{\Delta}_{ij} \subset \PP^1(\CC)^5$ is Zariski closed in $\PP^1(\CC)^5$, we see that $\Delta_{ij}$ is Zariski closed in $\ca N$. Define a Zariski open subset $\ca N_0  \subset \ca N$ as follows:
\begin{align*}
\ca N_0 \coloneqq \ca N - \bigcup_{i < j} \Delta_{ij}. 
\end{align*}

\begin{lemma} \label{lemma:coveringmarkingidentification}
The complex manifold $\ca N_0$ is connected and there is a canonical bijection between $\ca N_0$ and the set of marked smooth complex binary quintics $(F,m)$. Under this bijection, the map 
\begin{align} \label{align:GALOIS-cov}
\ca N_0 \to X_0(\CC), \quad \left((a_1, b_1), (a_2, b_2), \dotsc, (a_5, b_5) \right) \mapsto \prod_{j = 1}^5 \left(b_j\cdot x - a_j \cdot y \right)
\end{align}
corresponds to the map \eqref{markedquintics} that forgets the marking, and \eqref{align:GALOIS-cov} is a Galois covering with Galois group $\mf S_5$. 
\end{lemma}
\begin{proof}
Clear.
\end{proof}
In a similar way, we define, for $i,j,k \in \set{1, \dotsc, 5}$ with $i<j<k$, closed subsets
\[
\overline\Delta_{ijk} \coloneqq \set{(x_1, \dotsc, x_5) \in \PP^1(\CC)^5 \mid x_i = x_j = x_k}, \quad \quad \Delta_{ijk} \coloneqq f^{-1}(\overline\Delta_{ijk}) \subset \ca N,
\]
and put 
\[
\ca N_s \coloneqq \ca N - \bigcup_{i<j<k} \Delta_{ijk}. 
\]
The space $\ca N_s $ is equipped with a finite ramified covering map 
\[
\ca N_s \to X_s(\CC), \quad \quad \left((a_1, b_1), (a_2, b_2), \dotsc, (a_5, b_5) \right)\mapsto \prod_{j=1}^5 (b_j \cdot x- a_j \cdot y) \in \CC[x,y], 
\]
that commutes with \eqref{align:GALOIS-cov} and the natural open embedding $\ca N_0 \subset \ca N_s$. Observe that $\ca N_s \to X_s(\CC)$ is the Fox completion (cf.\ \cite{fox} or \cite[Section 8.1]{DeligneMostow}) of the spread $\ca N_0 \to X_0(\CC) \hookrightarrow X_s(\CC)$.

\subsection{Monodromy} 

Choose a marking $m_0$ lying over our base point $F_0 \in X_0(\CC)$.  
By Lemma \ref{lemma:coveringmarkingidentification}, 
we have a surjective homomorphism
\begin{equation}\label{eq:monodromy3}
\rho_5 \colon \pi_1(X_0(\CC), F_0) \twoheadrightarrow
\mf S_5
\end{equation}
whose kernel is given by the image of the natural embedding $\pi_1(\ca N_0, m_0) \hookrightarrow \pi_1(X_0(\CC), F_0)$. 
Composing the latter with the maps $\rho_\Gamma$ and $P(\rho_\Gamma)$ of \eqref{eq:monodromy-nonprojective} and \eqref{eq:monodromy} yields homomorphisms 
\begin{equation}\label{eq:monodromy2}
\mu\colon\pi_1(\ca N_0, m_0) \to \Gamma \quad\quad \text{and} \quad\quad P(\mu) \colon \pi_1(\ca N_0, m_0) \to P\Gamma.
\end{equation}
Consider the three-dimensional $\bb F_5$ vector space $\Lambda/(1 - \zeta)\Lambda$, as well as the quadratic space 
$$
W \coloneqq \left(\Lambda/(1-\zeta)\Lambda, \mf q \right)
$$
over $\bb F_5$. Here, $\mf q$ is the quadratic form obtained by reducing $\mf h$ modulo $(1-\zeta)\Lambda$. Define two groups $\Gamma_\theta$ and $P\Gamma_\theta$ as follows:
\[
\Gamma_\theta = \Ker\left( \Gamma \to \Aut(W) \right), \quad P\Gamma_\theta = \Ker\left( P\Gamma \to P\Aut(W) \right). 
\]
The following proposition seems to be due to Terada \cite{terada} and Yamazaki--Yoshida \cite{yamazakiyoshida}. 

\begin{proposition} 
 \label{prop:commutativemonodromy}
The image of the map $P(\mu)$ defined in (\ref{eq:monodromy2}) is the group $P\Gamma_\theta$. Moreover, the map $P\Gamma \to P\Aut(W)$ is surjective, and the induced homomorphism $$\overline{\rho}_\Gamma \colon \mf S_5 = \pi_1(X_0(\CC), F_0)/ \pi_1 \left( \ca N_0, m_0 \right)\to P\Gamma /P\Gamma_\theta = P\Aut(W)$$ is an isomorphism. As a consequence, we obtain a commutative diagram with exact rows:
\begin{equation} \label{eq:commutativemonodromy}
\begin{split}
    \xymatrix{
    0 \ar[r] & \pi_1(\ca N_0, m_0) \ar@{->>}[d]^{P(\mu)}\ar[r] & \pi_1(X_0(\CC), F_0)\ar[d]^{P(\rho_\Gamma)} \ar[r]^-{\rho_5} & \mf S_5\ar[d]_{\wr}^{\overline{\rho}_\Gamma} \ar[r] & 0 \\
    0 \ar[r] & P\Gamma_\theta \ar[r] & P\Gamma \ar[r] & P\Aut(W) \ar[r] & 0.
    }
    \end{split}
\end{equation}
\end{proposition}
\begin{proof}
The equality $P(\mu)\left(\pi_1(\ca N_0,m_0)\right) = P\Gamma_\theta$ follows from $\mu\left(\pi_1(\ca N_0,m_0)\right) = \Gamma_\theta$. For the latter, see 
\cite[Propositions 4.2 \& 4.3]{yamazakiyoshida}. The group $P\Aut(W) \cong \rm{PO}_3(\bb F_5)$ is isomorphic $\mf S_5$, and $\overline{\rho}_\Gamma \colon \mf S_5 \to P\Gamma/\Gamma_\theta$ is an isomorphism (cf.\ \cite[Propositions 4.2 \& 4.3]{yamazakiyoshida} and \cite[p.\ 10]{aperyyoshida-pentagonalstructure}). 
\end{proof}

\begin{corollary} \label{corollary:surjective-monodromy}
The monodromy representation $P(\rho_\Gamma) \colon \pi_1(X_0(\CC), F_0) \to P\Gamma$ is surjective. \hfill \qed
\end{corollary}



\subsection{Framed binary quintics, nodal binary quintics and local monodromy} 

By a \textit{framing} of a point $F \in X_0(\CC)$ we mean a projective equivalence class $[f]$, where \[
f \colon  \VV_F = \rm H^1(C_F(\CC), \ZZ) \xrightarrow{\sim} \Lambda\] is an $\OO_K$-linear isometry: two such isometries are in the same class if and only if they differ by an element in $\mu_K = (\OO_K^\ast)_{\rm{tors}}$. Let $\ca F_0$ be the collection of all framings of all points $F \in X_0(\CC)$. The set $\ca F_0$ can naturally be given the structure of a complex manifold, in a way similar to the procedure described in \cite[(3.9)]{ACTsurfaces}. In the sequel, we consider $\ca F_0$ as a complex manifold.
 
 \begin{lemma} \label{lemma:isomorphiccoveringspaces} 
\begin{enumerate}
 \item \label{item:connectedcovering-space}
The complex manifold $\ca F_0$ is connected and the map \begin{align} \label{align:F0}\ca F_0 \to X_0(\CC), \quad \quad (F, [f]) \mapsto F\end{align} that forgets the framing is a Galois covering map, with Galois group $P\Gamma$. 
\item \label{item:isomorphiccoveringspaces}
 The spaces $P\Gamma_\theta \setminus \ca F_0$ and $ \ca N_0$ are isomorphic as covering spaces of $X_0(\CC)$. In particular, there is a covering map $\ca F_0 \to \ca N_0$ with Galois group $P\Gamma_\theta$. 
\end{enumerate}
 \end{lemma}


\begin{proof}
Indeed, $\ca F_0$ is connected by Corollary \ref{corollary:surjective-monodromy}. The isomorphism $P\Gamma_\theta \setminus \ca F_0 \cong \ca N_0$ of covering spaces of $X_0(\CC)$ follows from the isomorphism $P\Gamma/P\Gamma_\theta \cong \mf S_5$ as quotients of $\pi_1(X_0(\CC), F_0)$, which was shown in Proposition \ref{prop:commutativemonodromy}. The lemma follows. 
\end{proof}



\begin{lemma} \label{lemma:irreduciblenormal}
The subvariety $\Delta \coloneqq X_s - X_0$ is an irreducible normal crossings divisor in $X_s$. 
\end{lemma}
\begin{proof}
The irreducibility of $\Delta$ is well-known, see e.g.~\cite[Section 14.1.1]{voisin}. 
Let $(p,F)$ be a point on the incidence variety $\ca I = \set{(p, F) \in \PP^1 \times X_s \mid p \in \tn{Sing}(F)}$, choose a hyperplane away from $p$, and view $p \in \CC$ and $F$ as a polynomial $f(x)$. 
Since $f^{''}(p)$ is non-zero, the tangent map $$Tq \colon T_{(p,f)} (\ca I) \to T_f(X_s) = X(\CC)$$ of the projection $q \colon \ca I \to X_s$ is injective, and its image consists of the linear subspace $X_p \subset X$ of binary quintics that contain $p$ as a root. If $f$ has $k$ double points $p_j$, where $k \in \set{1,2}$, then $\Delta$ is locally isomorphic to the union $\cup_{j=1}^k X_{p_j}$. These $X_{p_j}$ intersect transversally, and we are done. 
\end{proof}

\begin{definition} \label{def:deltak}
A \emph{node} $p \in Z_F(\CC) \subset \PP^1(\CC)$ of a binary quintic $F$ is a double point, i.e.\ a zero with multiplicity two. 
For $k = 1, 2$, let $\Delta_k \subset \Delta = X_s - X_0$ be the locus of stable binary quintics with exactly $k$ nodes. 
\end{definition}
The following result is due to Deligne and Mostow \cite{DeligneMostow}. 

\begin{lemma} \label{lemma:monodromy}
The local monodromy transformations of $\ca F_0 \to X_0(\CC)$ around every $F \in \Delta$ are of finite order. More precisely, if $F \in \Delta_k$, 
then the local monodromy group around $F$ is isomorphic to $(\ZZ/10)^k$. 
\end{lemma}
\begin{proof}
Let $F_1 \in \Delta_1$ be a binary quintic with one node. 
Consider the universal family $\mr C \to X_s$ of quintic covers of $\PP^1$ ramified along a stable binary quintic, whose restriction to $X_0$ is \eqref{universalcycliccover}. Let $D \subset X_s(\CC)$ be an open disc transverse to $\Delta_1$ at $F_1$. 
For $F \in D^\ast$, one obtains a monodromy transformation $T \colon \rm H^1(C_F(\CC), \ZZ) \to \rm H^1(C_F(\CC), \ZZ)$ induced by a vanishing cycle, and $T$ has order ten by \cite[Proposition 9.2]{DeligneMostow}. 
Similarly, if $F_2 \in \Delta_2$ has two nodes, we may choose an embedding $D^2 \subset X_s(\CC)$ of the polydisc $D^2 = D \times D$ transversal to $\Delta_2$ at $F_2$. Since distinct nodes have orthogonal vanishing cycles, the local monodromy transformations commute.
\end{proof}
In the following corollary, we let $D = \set{z \in \CC \colon \va{z} < 1}$ denote the open unit disc, and $D^\ast = D - \{0\}$ the punctured open unit disc. 

\begin{corollary} \label{cor:framedquintics}
There is an essentially unique branched cover 
\[
\pi \colon \ca F_s \to X_s(\CC),
\]
with $\ca F_s$ a complex manifold, 
such that for any $x \in \Delta$, any open $x \in U \subset X_s(\CC)$ with 
\[
U \cong D^k \times D^{6-k} \quad \text{and} \quad U \cap X_0(\CC) \cong (D^\ast)^k \times D^{6-k},\] 
and any component $U'$ of $\pi^{-1}(U) \subset \ca F_s$, there is an isomorphism $U' \cong D^k \times D^{6-k}$ such that the composition 
\[
D^k \times D^{6-k} \cong U' \to U \cong D^6 \; \text{ is the map } \; (z_1, \dotsc, z_6) \mapsto (z_1^{10}, \dotsc, z_k^{10}, z_{k+1}, \dotsc, z_6). 
\]
\end{corollary}
\begin{proof}
In view of \cite[Lemma 7.2]{beauvillecubicsurfaces} (see also \cite{fox} and \cite[Section 8.1]{DeligneMostow}), this follows from Lemma \ref{lemma:monodromy}. 
\end{proof}
The group $G(\CC) = \GL_2(\CC)/D$ (see equation (\ref{def:GC})) acts on $\ca F_0$ over its action on $X_0$. Explicitly, if $g \in G(\CC)$ and if $([\phi],  \phi\colon\VV_{F} \xrightarrow{\sim} \Lambda)$ is a framing of $F \in X_0(\CC)$, then 
\[
\left([\phi \circ g^\ast], \phi \circ g^\ast \colon \VV_{g\cdot F} \xrightarrow{\sim} \Lambda\right)
\] is a framing of $g\cdot F \in X_0(\CC)$. This is a left action. 
The group $P\Gamma$ also acts on $\ca F_0$ from the left, and the actions of $P\Gamma$ and $G(\CC)$ on $\ca F_0$ commute. By functoriality of the Fox completion, the action of $G(\CC)$ on $\ca F_0$ extends to an action of $G(\CC)$ on $\ca F_s$.

\begin{lemma} \label{lemma:freeaction}
The group $G(\CC) = \GL_2(\CC)/D$ acts freely on $\ca F_s$. 
\end{lemma}
\begin{proof}
Consider the natural action of $G(\CC)$ on $\ca N_s$, and the action of $G(\CC)$ on $P\Gamma_\theta \sm \ca F_s$ induced by the action of $G(\CC)$ on $\ca F_s$. With respect to these actions, 
the isomorphism of ramified covering spaces 
$P\Gamma_\theta \setminus \ca F_s \cong \ca N_s$ of $X_s(\CC)$ that results from Lemma \ref{lemma:isomorphiccoveringspaces}.\eqref{item:isomorphiccoveringspaces} is $G(\CC)$-equivariant. In particular, the natural ramified covering map $\ca F_s \to \ca N_s$ is $G(\CC)$-equivariant, 
and so it suffices to show that $G(\CC)$ acts freely on $\ca N_s$. 

To this end, note that $\ca N_s$ admits a natural $\CC^\ast$-quotient map \begin{align}\label{align:Cstarquotient}\ca N_s \to P_s,\end{align} where $P_s \subset \PP^1(\CC)^5$ is the space of stable ordered five-tuples in $\PP^1(\CC)$ introduced in Section \ref{realbinaryintroduction}, and where $\CC^\ast$ acts on $\ca N_s$ by $\lambda \cdot (v_1, v_2, \dotsc, v_5) = (\lambda v_1, v_2, \dotsc, v_5) =  (v_1, \lambda v_2, \dotsc, v_5) = \cdots = (v_1, v_2, \dotsc, \lambda v_5) \in \ca N_s$ for $\lambda \in \CC^\ast$ and $(v_1, v_2, \dotsc, v_5) \in \ca N_s$ (see the description of $\ca N_s$ in Section \ref{sec:marked}), and \eqref{align:Cstarquotient} is equivariant for the natural homomorphism $G(\CC) \to \PGL_2(\CC)$. 
Let $g \in \GL_2(\CC)$ and $x \in \ca N_s$ such that $gx = x$. Since any element of $ \PGL_2(\CC)$ that fixes three distinct points on $\PP^1(\CC)$ is the identity, we have that $\PGL_2(\CC)$ acts freely on $P_s$. Therefore, $g \in \CC^\ast \subset \GL_2(\CC)$. Let $F \in X_s(\CC)$ be the image of $x \in \ca N_s$; then $$F(x,y) = gF(x,y) = F(g^{-1}(x,y)) = F(g^{-1}x, g^{-1}y) = g^{-5}F(x,y).$$ Thus, we have $g^{5} = 1 \in \CC^\ast$, and we conclude that $g \in D \subset \GL_2(\CC)$. 
\end{proof}

\section{Deligne--Mostow uniformization of the moduli space of complex binary quintics} \label{complexball}

In this section, we show that results of 
Deligne and Mostow \cite{DeligneMostow} 
yield an isomorphism of complex analytic spaces $\ca M_s(\CC) \cong P\Gamma \setminus \CC H^2$. This map is induced by the Riemann extension $\ca P_s \colon \ca F_s \to \CCH^2$ of a holomorphic map $\ca P \colon \ca F_0 \to \CCH^2$ whose definition follows rather directly from the set-up and results of the previous section. 
We also show that this isomorphism induces an isomorphism between the divisor $G(\CC) \setminus \Delta(\CC) = \ca M_s(\CC) - \ca M_0(\CC)$ and the divisor $P\Gamma \setminus \mr H \subset P\Gamma \setminus \CC H^2$ defined by a certain hyperplane arrangement $\mr H \subset \CCH^2$, and prove that $\ca P_s$ identifies binary quintics with $k$ nodes ($k = 1,2$) with points in $\mr H$ where exactly $k$ hyperplanes meet. 
\subsection{The period map}  \label{section:periodmap}

Define a complex hermitian vector space $V$ as 
\[
V = \left( \Lambda \otimes_{\OO_K, \sigma_1} \CC, \mf h^{\sigma_1}\right),
\] 
where $\mf h^{\sigma_1}$ is the hermitian form defined in Section \ref{subsection:signature}. Let $\CC H^2$ be the space of negative lines in $V$. Using \cite[Proposition 11.7]{degaayfortman-nonarithmetic} and Lemma \ref{lemma:signature}, we see that the abelian scheme $J \to X_0$ induces a holomorphic map, the \emph{period map}:
\begin{equation} \label{eq:periodframed}
\ca P \colon \ca F_0 \to \CC H^2. 
\end{equation} Explicitly, if $(F, [f]) \in \ca F_0$ is the framing $[f\colon\rm H^1(C_F(\CC), \ZZ) \xrightarrow{\sim}\Lambda]$ of the binary quintic $F \in X_0(\CC)$, and if $JC_F$ is the Jacobian of the curve $C_F$, then 
$$f \left( \rm H^{0,-1}(JC_F)_{\zeta} \right) = f \left(\rm H^{1,0}(C_F)_{\zeta} \right) \subset \rm H^1(C(\CC), \CC)_{\zeta} = \Lambda \otimes_{\OO_K, \sigma_1} \CC = V$$
is a negative line in $V$, and we have $\ca P(F , [f] ) =  f \left(\rm H^{1,0}(C_F)_{\zeta} \right)  \in \CC H^2$. 
The map $\ca P$ is holomorphic, 
and descends to a morphism of complex analytic spaces
\begin{equation*} 
\ca M_0(\CC) =G(\CC) \setminus X_0(\CC) \to P\Gamma \setminus \CC H^2. 
\end{equation*}
Moreover, by Riemann extension, (\ref{eq:periodframed}) extends to a $G(\CC)$-equivariant holomorphic map
\begin{align}\label{eq:stableperiodmapframed}
\ca P_s\colon\ca F_s \to \CC H^2. 
\end{align}

\begin{theorem}[Deligne--Mostow] \label{th:delignemostow}
The period map (\ref{eq:stableperiodmapframed}) induces an isomorphism of complex manifolds
\begin{align} \label{eq:isomstablefivepoints-zero}
G(\CC) \setminus \ca F_s \cong \CC H^2. 
\end{align}
Taking $P\Gamma$-quotients gives an isomorphism of complex analytic spaces
\begin{equation} \label{eq:isomstablefivepoints}
\ca M_s(\CC) = G(\CC) \setminus X_s(\CC) \cong P\Gamma \setminus \CC H^2. 
\end{equation} 
\end{theorem}
\begin{proof}
Recall that $P_0  \subset \PP^1(\CC)^5$ is the set of $(x_1, \dotsc, x_5) \in \PP^1(\CC)^5$ such that all $x_i$ are distinct. In accordance with \cite{DeligneMostow}, define $Q = G(\CC) \setminus \ca N_0 = \PGL_2(\CC) \setminus P_0$ and $Q_{\textnormal{st}} \coloneqq G(\CC) \setminus \ca N_s = \PGL_2(\CC) \setminus P_s$. Fix a base point $0 \in Q$ whose image in $\mf S_5 \sm Q$ coincides with the image of $F_0 \in X_0(\CC)$ under the canonical map $X_0(\CC) \to \ca M_0(\CC) = \mf S_5 \sm Q$. 

By Lemma \ref{lemma:isomorphiccoveringspaces}, we have that $\ca F_0$ is a covering space of $\ca N_0$, with Galois group $P\Gamma_\theta$. 
In \cite{DeligneMostow}, Deligne and Mostow define $\widetilde Q \to Q$ to be the covering space corresponding to the monodromy representation $\pi_1(Q,0) \to P\Gamma$; since the image of this homomomorphism is $P\Gamma_\theta$ (see Proposition \ref{prop:commutativemonodromy}), it follows that $G(\CC) \setminus \ca F_0 \cong \widetilde Q$ as covering spaces of $Q$. Consequently, if $\widetilde Q_{\textnormal{st}} \to Q_{\rm{st}} $ denotes the Fox completion (cf.~\cite{fox}, \cite[Section 8.1]{DeligneMostow}) of the spread $$\widetilde Q \to Q \hookrightarrow  Q_{\textnormal{st}},$$ then there is an isomorphism $G(\CC) \setminus \ca F_s \cong \widetilde  Q_{\textnormal{st}}$ of branched covering spaces of $ Q_{\textnormal{st}}$. We obtain the following commutative diagram, in which the horizontal arrows on the left are isomorphisms:
$$
\xymatrixcolsep{5pc}
\xymatrix{
G(\CC) \setminus \ca F_s \ar[r]^-{\sim} \ar[d] & \widetilde{Q}_{\textnormal{st}}\ar[d] \ar[r] & \CC H^2  \ar[d] \\
G(\CC) \setminus \ca N_s \ar[r]^-{\sim}  \ar[d] & Q_{\textnormal{st}}  \ar[r] \ar[d]& P\Gamma_{\theta} \setminus \CC H^2\ar[d]\\
G(\CC) \setminus X_s(\CC) \ar[r]^-{\sim} &  \mf S_5 \setminus Q_{\textnormal{st}} \ar[r] & P\Gamma \setminus \CC H^2.
}
$$
The map $\widetilde{Q}_{\textnormal{st}} \to \CC H^2$ is an isomorphism by \cite[(3.11)]{DeligneMostow}. It follows that $Q_{\text{st}} \to P\Gamma_\theta \setminus \CC H^2$ and $\mf S_5 \setminus Q_{\text{st}} \to P\Gamma \setminus \CC H^2$ are isomorphisms as well. 
Therefore, we are done if we can show that the composition $\ca F_0 \to \widetilde Q \to \CC H^2$ agrees with the period map $\ca P$ of equation (\ref{eq:periodframed}). This follows from \cite[(2.23) and (12.9)]{DeligneMostow}. 
\end{proof}

\subsection{Nodal binary quintics and orthogonal hyperplanes}

Consider the CM type $\Phi = \set{\sigma_1, \sigma_2} \subset \Hom(K,\CC)$ defined in \eqref{equation:CMtype}, the hermitian $\OO_K$-lattice $(\Lambda, \mf h)$ defined in (\ref{eq:hermitianformonbinaryquinticlattice}), and the following sets (cf.~\cite[Sections 2.2 \& 2.3]{degaayfortman-nonarithmetic}): \begin{align} \label{align:caH}\ca H = \set{H_r \subset \CCH^2 \mid r \in \mr R}, \quad \tn{ and } \quad \mr H = \bigcup_{H\in \ca H}H \subset \CCH^2.\end{align} Here, $\mr R\subset \Lambda$ is the set of short roots, i.e.~the set of $r \in \Lambda$ with $\mf h(r,r) = 1$, and for each $r \in \mr R$, $H_r \subset \CC H^2$ is the hyperplane of elements $x \in \CC H^2$ that are orthogonal to $r$. 

\begin{proposition}\label{prop:conditions-unitary-orthogonal}
The hyperplane arrangement $\mr H \subset \CCH^2$ is an orthogonal arrangement in the sense of \cite{orthogonalarrangements}. In other words: any two different hyperplanes $H_1,H_2\in\ca H$ either meet orthogonally, or not at all. 
\end{proposition}
\begin{proof}
By Lemma \ref{lemma:signature}, we have that $\mf h$ has signature $(2,1)$ with respect to the embedding $\sigma_1 \colon K \hookrightarrow \CC$, and signature $(3,0)$ with respect to the embedding $\sigma_2 \colon K \hookrightarrow \CC$, where $\sigma_1$ and $\sigma_2$ are defined in \eqref{equation:CMtype}.  
Therefore, the result follows from \cite[Theorem 2.5]{degaayfortman-nonarithmetic} and \cite[Example 2.12]{degaayfortman-nonarithmetic}. 
\end{proof}

\begin{proposition} \label{prop:stableperiodshyperplane}
The map (\ref{eq:isomstablefivepoints}) induces an isomorphism of complex analytic spaces
\begin{equation*}
\ca M_0(\CC) = G(\CC) \setminus X_0(\CC) \cong P\Gamma \setminus \left(\CC H^2 - \mr H \right). 
\end{equation*}
\end{proposition}
\begin{proof}
We have $\ca P_s(\ca F_0) \subset \CC H^2 - \mr H$ by \cite[Proposition 11.12]{degaayfortman-nonarithmetic}, because the Jacobian of a smooth curve cannot contain a non-trivial abelian subvariety whose induced polarization is principal. Therefore, we have $\ca P_s^{-1}(\mr H) \subset \ca F_s - \ca F_0$. Since $\ca F_s$ is irreducible (it is smooth by Corollary \ref{cor:framedquintics} and connected by Lemma \ref{lemma:isomorphiccoveringspaces}.\eqref{item:connectedcovering-space}), the analytic space $\ca P_s^{-1}(\mr H)$ is a divisor. Since $\ca F_s - \ca F_0$ is also a divisor by Corollary \ref{cor:framedquintics}, we have $\ca P_s^{-1}(\mr H) = \ca F_s - \ca F_0$. 
\end{proof}
Define $\widetilde \Delta = \ca F_s - \ca F_0$, and for $k \in \set{1,2}$, let  $\widetilde \Delta_k \subset \widetilde \Delta$ be the inverse image of $\Delta_k$ in $\widetilde \Delta$ under the map $\widetilde \Delta \to \Delta$. Here, $\Delta_k \subset \Delta$ is the subvariety defined in Definition \ref{def:deltak}. 
Moreover, for $k = 1,2$, define $\mr H_k \subset \mr H$ as the set $$\mr H_k \coloneqq \set{x \in \CCH^2 \colon \va{\ca H(x)} = k}.$$ Thus, this is the locus of points in $\mr H$ where exactly $k$ hyperplanes meet. For $r \in \mr R$, define an isometry $$\phi_r \colon \Lambda \to \Lambda, \quad \quad \phi_r(x) = x - (1-\zeta) \mf h(x,r) \cdot r.$$ Let $[\phi_r] \in P\Gamma$ be the image of $\phi_r \in \Gamma$ in the group $P\Gamma = \Gamma/\mu_K$, and for $x \in \CCH^2$, define
\begin{align} \label{align:G(x)}
G(x) \coloneqq \langle [\phi_r] \mid r\in \mr R \mid \mf h(x,r) = 0 \rangle \subset P\Gamma. 
\end{align}
\begin{lemma} 
\label{lemma:stabilizergroups}
\begin{enumerate}
\item  \label{item:first-first}
The period map $\ca P_s$ of (\ref{eq:stableperiodmapframed}) satisfies $\ca P_s(\widetilde \Delta_k) \subset \mr H_k$. 
\item \label{item:seconc-seconc}
Let $\tilde F \in \widetilde \Delta_k \subset \ca F_s$ and $x = \ca P_s(\tilde F) \in \mr H_k \subset \CCH^2$. 
Let $P\Gamma_{\tilde F} \subset P\Gamma$ be the stabilizer of $\tilde F$ in $P\Gamma$. 
Then $P\Gamma_{\tilde F}= G(x)$, where $G(x) \cong (\ZZ/10)^k$ is as in \eqref{align:G(x)}. 
\end{enumerate}
\end{lemma}

\begin{proof}
\eqref{item:first-first}. We know that $\ca P_s$ induces an isomorphism $G(\CC) \setminus \widetilde \Delta \xrightarrow{\sim} \mr H$ by Theorem \ref{th:delignemostow} and Proposition \ref{prop:stableperiodshyperplane}. This map must identify the smooth (resp.~singular) locus of both analytic varieties with each other, from which the result follows.

\eqref{item:seconc-seconc}. This follows from Lemma \ref{lemma:monodromy} and Corollary \ref{cor:framedquintics} (compare \cite[(3.10)]{ACTsurfaces} and \cite[Lemma 10.3]{realACTsurfaces}). 
\end{proof}

\section{The moduli space of real binary quintics} \label{modrealbinquin}

With the period map for complex binary quintics in place, we turn to the construction of the period map for real binary quintics. 
Define $\kappa$ as the anti-holomorphic involution
\[
\kappa \colon X_0(\CC) \to X_0(\CC), \quad  F(x,y) = \sum_{i+j=  5}a_{ij}x^i y^j \;\; \mapsto \;\; \overline{F(x,y)}= \sum_{i+j =5}\overline{a_{ij}} x^i y^j. 
\]
\begin{definition}\label{definition:Gammaprime}  (Compare \cite[Sections 3.1 \& 3.2]{degaayfortman-nonarithmetic}.)
\begin{enumerate}
\item 
An $\OO_F$-linear bijection $\phi \colon \Lambda \xrightarrow{\sim} \Lambda$ is called \emph{anti-unitary} if $\phi(\mu \cdot x) = \overline \mu \cdot \phi(x)$ and $\mf h(\phi(x), \phi(y)) = \overline{\mf h(x,y)}$ for $x,y \in \Lambda$, $\mu \in \OO_K$. 
\item Let $\Gamma'$ be the group of unitary and anti-unitary $\OO_F$-linear bijections of $\Lambda$. Let $P\Gamma' = \Gamma'/\mu_K$. 
\item Let $\mr A$ be the set of anti-unitary involutions $\alpha \colon \Lambda \to \Lambda$, and define $P\mr A = \mu_K \setminus \mr A$. 
\item For $\alpha \in P\mr A$, define $\RRH^2_\alpha \subset \CCH^2$ as the fixed space of $\alpha$, i.e.\ $\RR H^2_\alpha = (\CC H^2)^\alpha$. 
\item For $\alpha \in P\mr A$, let $P\Gamma_\alpha \subset P\Gamma$ be the stabilizer of the subspace $\RRH^2_\alpha = (\CCH^2)^\alpha \subset \CCH^2$. 
\end{enumerate}
\end{definition}
For $\alpha \in P\mr A$, the notation $\RRH^2_\alpha$ reflects the fact that $\RRH^2_\alpha$ is isometric to the real hyperbolic plane $\RRH^2$, see \cite[Lemma 3.4]{degaayfortman-nonarithmetic}. 

For each $\alpha \in P\mr A$, there is a natural anti-holomorphic involution $\alpha \colon \ca F_0 \to \ca F_0$ lying over the anti-holomorphic involution $\kappa \colon X_0(\CC) \to X_0(\CC)$. 
To define $\alpha$, consider a framed binary quintic $(F, [f]) \in \ca F_0$, where $f\colon\VV_F \to \Lambda$ is an $\OO_K$-linear isometry. Let $C_F \to \PP^1_\CC$ be the induced quintic cover of $\PP^1_\CC$. Complex conjugation $\PP^1(\CC) \to \PP^1(\CC)$ induces induces a bijection $Z_F(\CC) \xrightarrow{\sim} Z_{\kappa(F)}(\CC)$, and extends to an anti-holomorphic diffeomorphism $
\sigma_F \colon C_F(\CC) \to C_{\kappa(F)}(\CC)$ 
with pull-back $\sigma_F^\ast\colon\bb V_{\kappa (F)} \to \bb V_F$. The composition $\alpha \circ f \circ \sigma_F^{\ast} \colon \bb V_{\kappa(F)} \to \Lambda$ 
induces a framing of $\kappa(F) \in X_0(\CC)$, and we define
\[
\alpha(F, [f]) \coloneqq \left(\kappa(F), [\alpha \circ f \circ \sigma_F^\ast] \right) \in \ca F_0.
\] Although we have chosen a representative $\alpha \in \mr A$ of the class $\alpha \in P\mr A$, the element $\alpha(F, [f]) \in \ca F_0$ does not depend on this choice.  

Consider the covering map $\ca F_0 \to X_0(\CC)$ defined in equation \eqref{align:F0} in Lemma \ref{lemma:isomorphiccoveringspaces}, and define $\ca F_0(\RR)$ as the preimage of $X_0(\RR)$ in the space $\ca F_0$. Then
\begin{align}\label{realpointsofFzero}
\ca F_0(\RR) = \coprod_{\alpha \in P\mr A} \ca F_0^\alpha \subset \ca F_0.
\end{align}
To see why the union in (\ref{realpointsofFzero}) is disjoint, let $\alpha \in P\mr A$; then 
$$
    \ca F_0^\alpha = \left\{ (F, [f]) \in \ca F_0 \colon\kappa(F) = F \textnormal{ and } [f \circ \sigma^\ast_{F} \circ f^{-1}] = \alpha \right\}.$$ Thus, for $\alpha, \beta \in P\mr A$ and $(F, [f]) \in \ca F_0^\alpha \cap \ca F_0^\beta$, we have $\alpha = [f \circ \sigma_F^\ast \circ f^{-1}] = \beta$. 

\begin{lemma} \label{lemma:silhol}
Let $\alpha \in P\mr A$. 
The anti-holomorphic involution $\alpha \colon \ca F_0 \to \ca F_0$ commutes with the period map $\ca P \colon \ca F_0 \to \CC H^2$ and the anti-holomorphic involution $\alpha \colon \CCH^2 \to \CCH^2$. 
\end{lemma}
\begin{proof}
If $\tn{conj} \colon \CC \to \CC$ is complex conjugation, then for any $F \in X_0(\CC)$, the induced map
$\sigma^\ast_F \otimes \tn{conj} \colon \bb V_{\kappa (F)}\otimes_\ZZ \CC \to \bb V_F \otimes_\ZZ \CC$ is anti-linear, preserves the Hodge decompositions \cite[Chapter I, Lemma 2.4]{silholsurfaces} as well as the eigenspace decompositions. 
\end{proof}
By Lemma \ref{lemma:silhol}, we obtain a \textit{real period map} 
\begin{align}\label{therealperiodmap}
\xymatrix{
\ca P^\RR \;\; \colon \;\; \ca F_0(\RR)  \ar@{=}[r] &\coprod_{\alpha \in P\mr A} \ca F_0^\alpha  \ar[r] &\coprod_{\alpha \in P\mr A} \RR H^2_\alpha  \ar@{=}[r]& 
\widetilde Y.
}
\end{align}
Let $\sigma \colon \GL_2(\CC) \to \GL_2(\CC)$ be the anti-holomorphic involution that sends a matrix to its complex conjugate, and note that $\sigma$ descends to an anti-holomorphic involution $\sigma \colon G(\CC) \to G(\CC)$; define $$G(\RR) \coloneqq G(\CC)^\sigma.$$
\begin{lemma}
The natural map $\GL_2(\RR) \to G(\RR)$ is an isomorphism. 
\end{lemma}
\begin{proof}
Indeed, if $M \in \GL_2(\CC)$ is a matrix such that $\sigma(M) = \zeta^j \cdot M$ for some $j \in \set{0, \dotsc, 4}$, then we can write $\zeta^j = \zeta^{2k}$ for some $k \in \set{0, \dotsc, 4}$ hence $\sigma(\zeta^k \cdot M) = \zeta^{-k}\zeta^j \cdot M = \zeta^k \cdot M$, so that $\zeta^k \cdot M \in \GL_2(\RR)$. This proves that $\GL_2(\RR) \to G(\RR)$ is surjective. Injectivity follows from the fact that the kernel of $\GL_2(\RR) \to G(\RR)$ is spanned by the elements of $D = \set{1, \zeta, \dotsc, \zeta^4}$ which are invariant under complex conjugation; the only such element is $1 \in D$.  
\end{proof}
The map (\ref{therealperiodmap}) is constant on $G(\RR)$-orbits, as the same is true for $\ca P \colon \ca F_0 \to \CC H^2$. By abuse of notation, we write $\RRH^2_\alpha - \mr H = \RRH^2_\alpha - \left(\mr H \cap \RRH^2_\alpha\right)$ for $\alpha \in P\mr A$.
\begin{proposition} \label{prop:realsmoothperiods} 
The period map (\ref{therealperiodmap}) descends to a $P\Gamma$-equivariant diffeomorphism 
\begin{align} \label{firstperiod}
\ca P^\RR \colon \ca M_0(\RR)^f \coloneqq G(\RR) \setminus \ca F_0(\RR) \cong \coprod_{\alpha \in P\mr A} \left(\RR H^2_\alpha - \mr H\right).
 \end{align}
 Let $C\mr A \subset P\mr A$ be a set of representatives for the action of $P\Gamma$ on $P\mr A$. 
By $P\Gamma$-equivariance, the map (\ref{firstperiod}) induces an isomorphism of real-analytic orbifolds
\begin{equation} \label{smoothrealperiodiso}
\ca P^\RR \colon \ca M_0(\RR) = G(\RR) \setminus X_0(\RR) \cong  \coprod_{\alpha \in C\mr A}P\Gamma_\alpha \setminus \left( \RR H^2_\alpha - \mr H
    \right). 
\end{equation}
\end{proposition}

\begin{proof}
We follow the proof of \cite[Theorem 3.3]{realACTsurfaces}. The first thing to observe is that the map $$\ca P^\RR \colon G(\RR) \setminus \ca F_0(\RR) \to \coprod_{\alpha \in P\mr A} \left(\RR H^2_\alpha - \mr H\right)$$ is a local diffeomorphism, as the same holds for $\ca P \colon G(\CC) \setminus \ca F_0 \to \CC H^2 - \mr H$ by Theorem \ref{th:delignemostow}. To prove the surjectivity and injectivity of $\ca P^\RR$, notice that the arguments used in \cite{realACTsurfaces} to prove the analogous claims for cubic surfaces readily carry over to our situation. 
\end{proof}
Our next goal is to prove the real analogue of Theorem \ref{th:delignemostow}. By the naturality of the Fox completion, for every $\alpha \in P\mr A$, the involution $\alpha \colon \ca F_0 \to \ca F_0$ extends to an involution $\alpha$ on $\ca F_s$. 
\begin{lemma} \label{lemma:alphaperiod} The restriction of the map $\ca P_s \colon \ca F_s \to \CC H^2$ to $\ca F_s^\alpha$ induces a diffeomorphism $$\ca P_s^{\alpha} \colon G(\RR) \setminus \ca F_s^\alpha \cong \RR H^2_\alpha.$$ 
\end{lemma}
\begin{proof}
The map $\ca P_s^{\alpha} \colon G(\RR) \setminus \ca F_s^\alpha \to \RR H^2_\alpha$ is a local diffeomorphism because its differential at any point is an isomorphism by Theorem \ref{th:delignemostow}. Let us prove that $\ca P_s^\alpha$ is injective. Apply \cite[Lemma 3.5]{realACTsurfaces} with $X = \ca F_0$, $G = G(\CC)$ and $\phi = \alpha$; then $X^\phi = \ca F_s^\alpha$ and $Z = G(\RR)$. Note that we may apply this lemma because $G(\CC)$ acts freely on $\ca F_s$, see Lemma \ref{lemma:freeaction}. The conclusion is that the map $$G(\RR) \setminus \ca F_0^\alpha \to G(\CC) \setminus \ca F_0 = \CC H^2 - \mr H$$ is injective. To prove the surjectivity of $\ca P_s^{\alpha}$, one uses \cite[Lemma 11.2]{realACTsurfaces} to see that the map $G(\RR) \setminus \ca F_s^{\alpha} \to G(\CC) \setminus \CC H^2$ is proper. By Proposition \ref{prop:realsmoothperiods}, its image contains the dense open subset $\ca P_s(\ca F_0^\alpha) = \RR H^2_\alpha - \mr H$, so $\ca P_s^\alpha$ is surjective.\end{proof}

\begin{definition}\label{definition:equivalence-relation}
(See \cite[Definition 1.1]{degaayfortman-nonarithmetic}.) Define an equivalence relation $\sim$ on the disjoint union $\coprod_{\alpha \in P\mr A}\RR H^2_\alpha$ in the following way. Consider two points $(x,\beta) \in \RRH^2_\beta \subset \coprod_{\alpha \in P\mr A}\RR H^2_\alpha$ and $(y,\gamma) \in \RRH^2_\gamma \subset \coprod_{\alpha \in P\mr A}\RR H^2_\alpha$. Then $(x,\beta) \sim (y,\gamma)$ if $x = y \in \CCH^2$ and $\gamma \circ \beta \in G(x)$.  
\end{definition}
By \cite[Lemma 4.4]{degaayfortman-nonarithmetic}, the action of $P\Gamma$ on $\coprod_{\alpha \in P\mr A}\RR H^2_\alpha$ is compatible with this equivalence relation; 
define 
$$Y \coloneqq \left(\coprod_{\alpha \in P\mr A}\RR H^2_\alpha\right) / \sim, \quad \quad M \coloneqq P\Gamma \sm Y.$$ 
\begin{theorem} \label{theorem:gluing-thm}
Let $C\mr A \subset P\mr A$ be a set of representatives for the action of $P\Gamma$ on $P\mr A$. There exists a canonical real hyperbolic orbifold structure on the topological space $M = P\Gamma \sm Y$ together with a natural open immersion of hyperbolic orbifolds 
\[
	\coprod_{\alpha \in C\mr A} P\Gamma_\alpha \sm \left( \RRH^2_\alpha - \mr H \right) \hookrightarrow M.
\]
Moreover, for each connected component $M_i \subset M$, there exists a lattice $P\Gamma_\RR^i \subset \rm{PO}(2,1)$ and an isomorphism of real hyperbolic orbifolds $M_i \cong P\Gamma_\RR^i \sm \RRH^2$. 
\end{theorem}
\begin{proof}
By Lemma \ref{lemma:signature} \& \cite[Example 2.12]{degaayfortman-nonarithmetic}, this is a special case of \cite[Theorem 1.2]{degaayfortman-nonarithmetic}. 
\end{proof}

Let $p \colon \coprod_{\alpha \in P\mr A}\RR H^2_\alpha \to Y$ be the quotient map, consider the map $\pi\colon\ca F_s \to X_s(\CC)$ (see Corollary \ref{cor:framedquintics}) and define a union of embedded real submanifolds of $\ca F_s$ as follows:
\[
\ca F_s(\RR) = \bigcup_{\alpha \in P\mr A} \ca F_s^\alpha = \pi^{-1}\left(X_s(\RR)\right).
\]
We arrive at the main theorem of Section \ref{modrealbinquin}. 
\begin{theorem} \label{th:realstableperiod}
There is a smooth map 
\begin{align}\label{therealstableperiodmap}
\ca P_s^{\RR} \colon  \coprod_{\alpha \in P\mr A} \ca F_s^\alpha \to  \coprod_{\alpha \in P\mr A}\RR H^2_\alpha 
\end{align} that extends the real period map (\ref{therealperiodmap}). The map (\ref{therealstableperiodmap}) induces the following commutative diagram of topological spaces, in which $\mr P_s^\RR$ and $\mr T^\RR_s$ are homeomorphisms:
\begin{equation*}
\xymatrixcolsep{5pc}
\xymatrix{
&\coprod_{\alpha \in P\mr A} \ca F_s^\alpha \ar[r]^{\ca P_s^{\RR}}\ar[d] &\coprod_{\alpha \in P\mr A} \RR H^2_\alpha \ar[d]^p\\
&\ca F_s(\RR) \ar[r]^{\ca P_s^\RR}\ar[d] & Y\ar@{=}[d] \\
\ca M_s(\RR)^f \ar@{=}[r]\ar[d] &G(\RR) \setminus \ca F_s(\RR) \ar[r]^{\mr P_s^\RR}_\sim\ar[d] & Y\ar[d] \\
\ca M_s(\RR) \ar@{=}[r] &G(\RR) \setminus X_s(\RR) \ar[r]^{{\mr T_s}^\RR}_\sim & P\Gamma \setminus Y. 
}
\end{equation*}
%
\end{theorem}

\begin{proof}
The existence of $\ca P_s^{\RR}$ follows from the compatibility between $\ca P_s$ and the involutions $\alpha \in P\mr A$. We claim that the composition $p \circ \ca P_s^{\RR}$ 
factors through $\ca F_s(\RR)$. To prove this, let $(f,\alpha)$ and $(g, \beta)$ be elements of the disjoint union $\coprod_{\gamma \in P\mr A} \ca F_s^\gamma$, with $f \in \RRH^2_\alpha$ and $g \in \RRH^2_\beta$. Then $(f, \alpha)$ and $(g,\beta)$ have the same image in $\ca F_s(\RR)$ if and only if $f = g \in \ca F_s^\alpha \cap \ca F_s^\beta$, in which case 
$
\ca P_s(f) = \ca P_s(g) \eqqcolon x \in \RR H^2_\alpha \cap \RRH^2_\beta$. Let $(x, \alpha)$ and $(y,\beta)$ be elements of the disjoint union $\widetilde Y = \coprod_{\gamma \in P\mr A}\RRH^2_\gamma$, with $x \in \RRH^2_\alpha$ and $y =x \in \RRH^2_\beta$. We need to prove $(x,\alpha) \sim (x,\beta) \in \widetilde Y$, for the equivalence relation $\sim$ on $\widetilde Y$ defined in Definition \ref{definition:equivalence-relation}. Note that $\alpha\beta \in P\Gamma_f$, 
and that $\ca P_s$ induces an isomorphism 
$
P\Gamma_f \cong G(x)$, see Lemma \ref{lemma:stabilizergroups}. Hence $\alpha \beta \in G(x)$ so that $(x,\alpha) \sim (x,\beta)$, proving what we want. We conclude that the composition $p \circ \ca P_s^{\RR}$ factors through a map $\ca P_s^\RR \colon \ca F_s(\RR) \to Y$. 

Next, we prove the $G(\RR)$-equivariance of $\ca P_s^{\RR}$. Suppose that 
$
f \in \ca F_s^\alpha, g \in \ca F_s^\beta$ such that $ a \cdot f = g \in \ca F_s(\RR)$ for some $ a \in G(\RR)$. Then $x\coloneqq \ca P_s(f) = \ca P_s(g) \in \CC H^2$, so we need to show that $\alpha \beta \in G(x)$. The actions of $G(\CC)$ and $P\Gamma$ on $\CC H^2$ commute, and the same holds for the actions of $G(\RR)$ and $P\Gamma'$ on $\ca F_s^\RR$, where $P\Gamma'$ is the group defined in Definition \ref{definition:Gammaprime}. 
It follows that 
$
\alpha(g) = \alpha (a \cdot f) = a \cdot \alpha(f) = a \cdot f = g,
$
hence $g \in \ca F_s^\alpha \cap \ca F_s^\beta$. This implies in turn that $(\alpha \circ \beta) (g) = g$, hence $\alpha \beta \in P\Gamma_g \cong G(x)$, see Lemma \ref{lemma:stabilizergroups}. Therefore, $(x,\alpha) \sim (x,\beta)$, so that $\alpha \beta \in G(x)$ as desired. 

Let us prove that $\mr P_s^\RR$ is injective. To do so, let 
$
f \in \ca F_s^\alpha$ and $g \in \ca F_s^\beta$ and suppose that these elements have the same image in $Y$. Thus, $
x\coloneqq \ca P_s(f) = \ca P_s(g) \in \RR H^2_\alpha \cap \RRH^2_\beta,
$
and $\beta = \phi \circ \alpha$ for some $\phi \in G(x)$. 
We have $\phi \in G(x) \cong P\Gamma_f$ (Lemma \ref{lemma:stabilizergroups}) hence $    \beta(f) = \phi \left(\alpha (f)\right) = \phi(f) = f$. 
Therefore $f,g \in \ca F_s^\beta$; since $\ca P_s(f) = \ca P_s(g)$, it follows from Lemma \ref{lemma:alphaperiod} that there exists $a \in G(\RR)$ such that $a \cdot f = g$. This proves the injectivity of $\mr P_s^\RR$.   
The surjectivity of $\mr P_s^\RR\colon G(\RR) \setminus \ca F_s(\RR) \to Y$ is straightforward: it follows from the surjectivity of $\ca P_s^{\RR}$, see Lemma \ref{lemma:alphaperiod}. Finally, we claim that $\mr P_s^\RR$ is open. Let $U \subset G(\RR) \setminus \ca F_s^\RR$ be open. 
Let $V$ be the preimage of $U$ in $\coprod_{\alpha \in P\mr A}\ca F_s^\alpha$. Then 
$
V = (\ca P_s^{\RR})^{-1}\left( p^{-1}\left(\mr P_s^\RR(U)\right)\right)
$, hence 
$
    \ca P_s^{\RR}\left( V \right) = p^{-1}\left(\mr P_s^\RR(U)\right)$, so that it suffices to show that $\ca P_s^\RR(V)$ is an open subset of $\coprod_{\alpha \in P\mr A} \RR H^2_\alpha$. This follows, because $\ca P_s^{\RR}$ is open, being the coproduct of the maps $\ca F_s^\alpha \to \RR H^2_\alpha$, which are open since they have surjective differential at each point. 
\end{proof}

\begin{corollary} \label{cor:theorem2}
 Let $C\mr A \subset P\mr A$ be a set of representatives for the action of $P\Gamma$ on $P\mr A$. Then there is a lattice $P\Gamma_\RR \subset \textnormal{PO}(2,1)$, an open immersion of hyperbolic orbifolds
 \begin{align}\label{inclusionofforbifodls}
 \coprod_{\alpha \in C\mr A}P\Gamma_\alpha \setminus \left( \RR H^2_\alpha - \mr H
    \right) \hookrightarrow P\Gamma_\RR \setminus \RR H^2, 
 \end{align}
 and a homeomorphism
\begin{equation} \label{stablehom}
\ca P_s^{\RR} \colon \ca M_s(\RR) =  G(\RR) \setminus X_s(\RR) \cong P\Gamma_\RR \setminus \RR H^2,
\end{equation}
such that the restriction of \eqref{stablehom} to $\ca M_0(\RR) \subset \ca M_s(\RR)$ coincides with the isomorphism (\ref{smoothrealperiodiso}). 
\end{corollary}
\begin{proof}
This follows directly from Theorems \ref{theorem:gluing-thm} and \ref{th:realstableperiod}.
\end{proof}

\begin{remark}
The proof of Corollary \ref{cor:theorem2} also shows that $\ca M_s(\RR)$ is homeomorphic to a complete hyperbolic orbifold 
in the cases where $\ca M_s$ is the stack of cubic surfaces or of binary sextics over $\RR$. This strategy to uniformize the real moduli space differs from the one used in \cite{realACTnonarithmetic, realACTbinarysextics,realACTsurfaces}, since we first glue the real ball quotients together (by using the general construction of \cite{degaayfortman-nonarithmetic}) and then prove that the real moduli space is homeomorphic to the resulting glued space.
\end{remark}


\section{The moduli space of real binary quintics as a hyperbolic triangle} \label{space:hyperbolictriangle}

Consider the moduli space $\ca M_s(\RR) = \GL_2(\RR) \setminus X_s(\RR)$ of stable real binary quintics. Let $\va{\ca M_s(\RR)}$ be the underlying topological space of $\ca M_s(\RR)$. 

\begin{definition} \label{neworbifold}
Let $\overline{\mr M}_\RR$ be the orbifold with $\va{\ca M_s(\RR)}$ as underlying space, whose orbifold structure is induced by the homeomorphism (\ref{stablehom}) and the natural orbifold structure of $P\Gamma_\RR \setminus \RRH^2$. 
\end{definition}
The goal of Section \ref{space:hyperbolictriangle} is to prove the following result. 

\begin{theorem}\label{theorem:lattice-thm}
Consider the lattice $P\Gamma_\RR \subset \rm{PO}(2,1)$, see Corollary \ref{cor:theorem2}, and the hyperbolic orbifold $\overline{\mr M}_\RR \cong P\Gamma_\RR \sm \RRH^2$, see Definition \ref{neworbifold}. Then $\overline{\mr M}_\RR$ is isometric to the hyperbolic triangle of angles $\pi/3, \pi/5, \pi/10$. In particular, $P\Gamma_\RR$ is conjugate to the lattice $\Gamma_{3,5,10}$ defined in \eqref{PGAMMAR}. 
\end{theorem}
To prove Theorem \ref{theorem:lattice-thm}, we need to 
understand the orbifold structure of $\ca M_s(\RR)$ and how this structure differs from the orbifold structure of the quotient space $P\Gamma_\RR \setminus \RR H^2$, see Corollary \ref{cor:theorem2}. To this end, we will first analyze the orbifold structure of $\ca M_s(\RR)$ by listing its stabilizer groups. 

\subsection{Automorphism groups of stable real binary quintics} \label{sec-sec:1} 

Recall that there is a canonical orbifold isomorphism
$$
\ca M_s(\RR) = G(\RR) \setminus X_s(\RR) = \PGL_2(\RR) \setminus (P_s/\mf S_5)(\RR).
$$
Thus, to list those groups that occur as the automorphism group of a binary quintic is to classify the stabilizer groups $\PGL_2(\RR)_x$ of points $x = (x_1, \dotsc, x_5) \in (P_s/\mf S_5)(\RR)$. 

\begin{proposition} \label{prop:complex-stabilizer}
Consider the stabilizer group $\PGL_2(\CC)_x$ of a point $x \in P_0/\mf S_5$. If $\PGL_2(\CC)_x$ is non-trivial, then $\PGL_2(\CC)_x$ is isomorphic to either $\ZZ/2, \ZZ/4, D_3$ or $D_5$. Moreover, the conjugacy class of each such a subgroup of $\PGL_2(\CC)$ is unique. If $H$ equals any of the subgroups $\ZZ/4$, $D_3$ or $D_5$ of $ \PGL_2(\CC)$, then there is a unique $\PGL_2(\CC)$-orbit in $P_0/\mf S_5$ with stabilizer group $H$. 
\end{proposition}
\begin{proof}

By \cite[Theorem 4.2]{beauvillePGL2}, any finite subgroup of $\PGL_2(\CC)$ is isomorphic to $\ZZ/n$, $D_n$ (the dihedral group of order $2n$), $\mf A_4$, $\mf S_4$ or $\mf A_5$, and there is only one conjugacy class for each of these groups. 
Let $H$ be any of these groups, considered as a subgroup of $\PGL_2(\CC)$. Assume that, with respect to the action of $H$ on the finite subsets of $\PP^1(\CC)$, one has
\[
H \cdot \set{z_1, \dotsc, z_5} = \set{z_1, \dotsc, z_5} \subset \PP^1(\CC), \quad \text{for} \quad z_1, \dotsc, z_5 \in \PP^1(\CC) \quad \text{distinct.}
\]
This gives a homomorphism $\rho \colon H \to \mf S_5$ as follows: for an element $j \in \set{1,2,3,4,5}$, we let $\rho(h)(j) \in \set{1,2,3,4,5}$ be the element with $z_{\rho(h)(j)} = h \cdot z_j$. 

Note that $\rho$ is injective, as $h \cdot z_i = z_i$ for each $i$ implies $h = \id$. Therefore, 
\[
H \in \set{\ZZ/2, \ZZ/3, \ZZ/4, \ZZ/5, \ZZ/6, D_3, D_4, D_5, \mf A_4, \mf S_4, \mf A_5}. 
\]
Next, assume that $H = \Stab_{\PGL_2(\CC)}(x)$ for the five-element subset $x = \set{z_1, \dotsc, z_5} \subset \PP^1(\CC)$. Suppose that $\phi \in H$ is an element of order three. Note that there must be three distinct elements $z_i \in x$ with $\phi(z_i) \neq z_i$. We may assume that these are $z_1, z_2$ and $z_3$. Moreover, we may assume that $\phi(z_1) = z_2$, $\phi(z_2) = z_3$ and $\phi(z_3) = z_1$. By replacing $\phi$ by $g \phi g^{-1}$ for some $g \in \PGL_2(\CC)$, we may assume that $z_1= 1, z_2 = \zeta_3$ and $z_3 = \zeta_3^2$, and that $\phi(z) = \zeta_3 \cdot z$ for $z \in \PP^1(\CC)$. This gives $x = \set{1, \zeta_3, \zeta_3^2, z_4, z_5} \subset \PP^1(\CC)$. As $\phi(z_4) \neq z_5$, we have $\phi(z_4) = z_4$ and $\phi(z_5) = z_5$, so that 
 \begin{align}\label{x}x = \set{1, \zeta_3, \zeta_3^2, 0, \infty} \subset \PP^1(\CC).\end{align} Let $\nu \in \PGL_2(\CC)$ be the element with $\nu(z) = 1/z$ for $z \in \PP^1(\CC)$. Then \eqref{x} implies that $\nu$ is contained in $H$, and one readily observes that $H = D_3$. We conclude that 
 \[
H = \text{Stab}_{\PGL_2(\CC)}(x) \in \set{\ZZ/2, \ZZ/4, \ZZ/5, D_3, D_4, D_5}. 
\]
It remains to exclude $\ZZ/5$ and $D_4$. Suppose that $H$ contains an element $\phi$ of order five. Similar to the above, one readily shows that one may assume that
$$
x = \set{1, \zeta_5, \zeta_5^2, \zeta_5^3, \zeta_5^4} \subset \PP^1(\CC), \quad \phi(z) = \zeta_5 \cdot z \quad \text{for} \quad z \in \PP^1(\CC).
 $$
This implies that the element $\nu \in \PGL_2(\CC)$ as defined above is contained in $H$, and $H = D_5$. Finally, assume $H$ contains an element $\phi$ of order four. We may assume that 
\[
x = \set{1, i, -1, -i, 0}, \quad \phi(z) = i \cdot z \quad \text{for} \quad z \in \PP^1(\CC). 
\]
As a consequence, we have $H = \ZZ/4$, and the proof is finished. 
\end{proof}
We proceed to prove the real analogue of Proposition \ref{prop:complex-stabilizer}.

\begin{proposition} \label{prop:allstabilizergroups}
Let $x \in (P_s/ \mf S_5)(\RR)$ such that the stabilizer group $\PGL_2(\RR)_x \subset \PGL_2(\RR)$ of $x$ is non-trivial. Then its stabilizer group $\PGL_2(\RR)_x$ is isomorphic to $\ZZ/2$, $D_3$ or $D_5$. For each $n \in \{3,5\}$, there is a unique $\PGL_2(\RR)$-orbit of points $x$ in $(P_s/ \mf S_5)(\RR)$ with stabilizer $D_n$.  
\end{proposition}
\begin{proof}
We have an injection 
$
(P_s/\mf S_5)(\RR) \hookrightarrow P_s/\mf S_5$ which is equivariant for the embedding $\PGL_2(\RR) \hookrightarrow \PGL_2(\CC)$. In particular, $\PGL_2(\RR)_x \subset \PGL_2(\CC)_x$ for every $x \in (P_s/\mf S_5)(\RR)$. 
Note that none of the groups appearing in Proposition \ref{prop:complex-stabilizer} have subgroups isomorphic to $D_2 = \ZZ/2 \rtimes \ZZ/2$ or $D_4 = \ZZ/2 \rtimes \ZZ/4$. Consider the involution $    \nu = (z \mapsto 1/z) \in \PGL_2(\RR)$. 
We will prove the proposition by using the following steps.   


\begin{enumerate}[wide, labelwidth=!, labelindent=0pt]
\item[\bf{Step 1:}] \label{item:involution}
Let $\tau \in \PGL_2(\RR)$. Consider a subset $S = \{x,y,z\} \subset \PP^1(\CC)$ stabilized by complex conjugation, such that $\tau(x) = x$, $\tau(y) = z$ and $\tau(z) = y$. There is a transformation $g \in \PGL_2(\RR)$ that maps $S$ to either $\{-1, 0, \infty\}$ or $\{-1, i, -i\}$, and that satisfies $g \tau g^{-1} = \nu = (z \mapsto 1/z) \in \PGL_2(\RR)$. In particular, $\tau^2 = \id$. 
\begin{proof}
Indeed, two transformations $g,h \in \PGL_2(\CC)$ that satisfy $g(x_i) = h(x_i)$ for three different points $x_1, x_2, x_3 \in \PP^1(\CC)$ are necessarily equal. 
\end{proof}

\item[\bf{Step 2:}] \label{item-two}
There is no $\phi \in \PGL_2(\RR)$ of order four that fixes a point $x \in (P_s/\mf S_5)(\RR)$. 
\begin{proof}
By \cite[Theorem 4.2]{beauvillePGL2}, all subgroups $G \subset \PGL_2(\RR)$ that are isomorphic to $\ZZ/4$ are conjugate to each other. Since the transformation $I\colon z \mapsto (z-1)/(z+1)$ is of order four, it gives a representative $G_I = \langle I \rangle$ of this conjugacy class. 
It is easily shown that $I$ cannot fix any point $x \in (P_s/\mf S_5)(\RR)$. 
\end{proof}
%

\item[\bf{Step 3:}] \label{item:D3}
Define
$\rho \in \PGL_2(\RR)$ by  $\rho(z) = \frac{-1}{z+1}$. 
Let $x = (x_1, \dotsc, x_5) \in (P_s/\mf S_5)(\RR)$. Let $\phi \in \PGL_2(\RR)$ of order three, with $\phi(x) = x$. There is a transformation $g \in \PGL_2(\RR)$ mapping $x$ to $(-1, \infty, 0, \omega, \omega^2)$ with $\omega$ a primitive third root of unity. The stabilizer of $x$ to the subgroup of $\PGL_2(\RR)$ generated by $\rho$ and $\nu$. In particular, we have $\PGL_2(\RR)_x \cong D_3$. 

\begin{proof}
By Step 1, there are elements $x_1, x_2, x_3$ which form an orbit under $\phi$. Since complex conjugation preserves this orbit, one element in it is real; since $g$ is defined over $\RR$, they are all real. Let $g \in \PGL_2(\RR)$ such that $g(x_1) = -1$, $g(x_2) = \infty$ and $g(x_3) = 0$. Define $\kappa = g \phi g^{-1}$. Then $\kappa^3 = \id$, and $\kappa$ preserves $\{-1, \infty, 0\}$ and sends $-1$ to $\infty$ and $\infty$ to $0$. Consequently, $\kappa(0) = -1$, and it follows that $\kappa = \rho$. Hence $x$ is equivalent to an element of the form 
$(-1, \infty, 0, \alpha, \beta)$ with $\beta = \bar \alpha$ and $\alpha^2 + \alpha + 1 = 0$. 
\end{proof}
Recall that $ \lambda = \zeta_5 + \zeta_5^{-1} \in \RR$. Define 
$
\gamma \in \PGL_2(\RR)$ by $\gamma(z) = \frac{(\lambda + 1)z - 1}{z + 1}$ for $ z \in \PP^1(\CC).
$

\item[\bf{Step 4:}] \label{item:D5}
Let $x = (x_1, \dotsc, x_5) \in (P_s/\mf S_5)(\RR)$. Suppose $x$ is stabilized by a subgroup of $\PGL_2(\RR)$ of order five. There is a transformation $g \in \PGL_2(\RR)$ mapping $x$ to $z = (0, -1, \infty, \lambda+1, \lambda)$ and identifying the stabilizer of $x$ with the subgroup of $\PGL_2(\RR)$ generated by $\gamma$ and $\nu$. In particular, the stabilizer $\PGL_2(\RR)_x$ of $x$ is isomorphic to $D_5$. 

\begin{proof}
Let $\phi\in \PGL_2(\RR)_x$ be an element of order five. 
Using Step 1, one shows that the $x_i$ are pairwise distinct, and we may assume that $x_i = \phi^{i-1}(x_1)$ for $i = 2, \dotsc, 5$. Since there is one real $x_i$ and $\phi$ is defined over $\RR$, all $x_i$ are real. 

Let $z = \set{0, -1, \infty, \lambda + 1, \lambda}$. Note that $z$ is the orbit of $0$ under $\gamma\colon z \mapsto ((\lambda+1)z - 1)/(z+1)$. The reflection $\nu\colon z \mapsto 1/z$ preserves $z$ as well: we have $\lambda + 1 = - (\zeta_5^2 + \zeta_5^{-2}) = - \lambda^2 + 2$, so that $\lambda(\lambda + 1) = 1$. We conclude that $\PGL_2(\RR)_z \cong D_5$. Thus, by Proposition \ref{prop:complex-stabilizer}, 
there exists $g \in \PGL_2(\CC)$ such that $g(x_1) = 0$, $g(x_2) = -1$, $g(x_3) = \infty$, $g(x_4) = \lambda + 1$ and $g(x_5) = \lambda$, and such that $g\PGL_2(\RR)_xg^{-1} = \PGL_2(\RR)_z$. Since all $x_i$ and $z_i \in z$ are real, we have $\bar g(x_i) = z_i$ for each $i$. Hence, $g$ and $\bar g$ coincide on three points, which implies that $g = \bar g$, i.e.\ $g  \in \PGL_2(\RR)$. 
\end{proof}
\end{enumerate}
By Steps 1--4 above, together with Proposition \ref{prop:complex-stabilizer}, we are done. 
\end{proof}


 \subsection{Comparing the orbifold structures} \label{sec-sec:3} 

%

There are two orbifold structures on the space $\va{\ca M_s(\RR)}$: on the one hand, one has the natural orbifold structure on $\ca M_s(\RR)$ by considering it as the real locus of a smooth separated Deligne--Mumford stack over $\RR$ (see \cite[Proposition 2.12]{thesis-degaayfortman}), this is the orbifold structure of the quotient $ G(\RR) \setminus X_s(\RR)$; on the other hand, one has the orbifold structure $\overline{\mr M}_\RR$ introduced in Definition \ref{neworbifold}. The goal of Section \ref{sec-sec:3} is to calculate the difference between these orbifold structures. 

We first show that there are no cone points in the orbifold $\PGL_2(\RR) \setminus (P_s/ \mf S_5)(\RR)$. These are orbifold points whose stabilizer group is $\ZZ/n$ ($n \geq 2$) acting on the orbifold chart by rotations. By Proposition \ref{prop:allstabilizergroups}, the fact that $\PGL_2(\RR) \setminus (P_s/ \mf S_5)(\RR)$ has no cone points follows from:
\begin{lemma} \label{lemma:zmod2stabilizer}
Let $x=(x_1, \dotsc, x_5) \in (P_s/\mf S_5)(\RR)$ such that $\PGL_2(\RR)_x = \langle \tau \rangle $ has order two. There is a $\PGL_2(\RR)_x$-stable open neighborhood $U \subset (P_s/\mf S_5)(\RR)$ of $x$ such that $\PGL_2(\RR)_x \setminus U \to \ca M_s(\RR)$ is injective, and a homeomorphism $\phi\colon(U,x) \to (B,0)$ for $0 \in B \subset \RR^2$ an open ball, such that $\phi$ identifies $\PGL_2(\RR)_x$ with $\ZZ/2$ acting on $B$ by reflections in a line through $0$.  
\end{lemma}
\begin{proof}
Using Step 1 in the proof of Proposition \ref{prop:allstabilizergroups}, one checks that the only possibilities for the element $x=(x_1, \dotsc, x_5) \in (P_s/\mf S_5)(\RR)$ are  $(-1, 0, \infty, \beta, \beta^{-1})$,  $(-1, i, -i, \beta, \beta^{-1}) $, $ (-1, -1, \beta, 0, \infty)$, $ (-1, -1, \beta, i, -i) $, $(0,0, \infty, \infty, -1)$ and $(-1, i, i, - i, -i)$ for some $\beta \in \PP^1(\RR)$. 
\end{proof}

To analyze the difference between the two orbifolds $\ca M_s(\RR)$ and $\overline{\mr M}_\RR$, we also need the following general lemma. 

\begin{lemma} \label{lemma:generallemma-groups}
Let $X$ be a set and let $\Gamma$ and $G$ be groups with commuting actions on $X$. Let $x \in X$ with images $\bar x \in \Gamma \sm X$ and $[x] \in G \sm X$. Let $\Gamma_{[x]}$ be the stabilizer of $[x] \in G \sm X$ in $\Gamma$, and let $G_{\bar x}$ be the stabilizer of $\bar x$ in $G$. Then for each $\gamma \in \Gamma_{[x]}$ there exists an element $\phi(\gamma) \in G_{\bar x}$, unique up to multiplication by an element of $G_x$, such that $\gamma \cdot x = \phi(\gamma) \cdot x$;  moreover, the map 
\begin{align} \label{align:group-set-lemma}
\Gamma_{[x]}/\Gamma_x \to G_{\bar x}/G_x, \quad \quad \gamma \mapsto \phi(\gamma)
\end{align}
is an isomorphism. 
\end{lemma}
\begin{proof}
The map \eqref{align:group-set-lemma} is well-defined because if $g, g' \in G$ are such that $\gamma \cdot x= g \cdot x = g' \cdot x$ then $(g')^{-1}g \in G_x$. Since the  construction is symmetric in $\Gamma$ and $G$, the analogous map $G_{\bar x}/G_x \to  \Gamma_{[x]}/\Gamma_x$ is also well-defined. The latter is a left and right inverse of \eqref{align:group-set-lemma}, and we are done.
\end{proof}

Recall, see \cite[Proposition 13.3.1]{Thurston80}, that the singular locus of a two-dimensional orbifold has the following types of local models: (i) $\RR^2/(\ZZ/2)$, where $\ZZ/2$ acts on $\RR^2$ by reflection in the $y$-axis (\emph{mirror points}); (ii) $\RR^2/(\ZZ/n)$, with $\ZZ/n$ acting by rotations (\emph{cone points of order $n$}); (iii) $\RR^2/D_n$, with $D_n$ the dihedral group of order $2n$ (\emph{corner reflectors of order $n$}). 

\begin{proposition} \label{prop:conesreflectors} Consider the orbifold structures $\ca M_s(\RR)$ and $\overline{\mr M}_\RR$ on the space $\va{\ca M_s(\RR)}$. 
\begin{enumerate}
\item For $x_0 \in \ca M_s(\RR)$, the isomorphism class of stabilizer groups of $\ca M_s(\RR)$ and $\overline{\mr M}_\RR$ at $x_0$ differ if and only if 
$x_0 \in \ca M_s(\RR)$ is the moduli point attached to the five-tuple $(\infty, i, i, -i, -i)$. 
\item The stabilizer group of $\ca M_s(\RR)$ at the point $x_0$ is isomorphic to $\ZZ/2$, whereas the stabilizer group of $\overline{\mr M}_\RR$ at $x_0$ is isomorphic to the dihedral group $D_{10}$ of order $20$. 
\item The orbifold $\overline{\mr M}_\RR$ has no cone points and three corner reflectors, whose angles are $\pi/3, \pi/5$ and $\pi/10$. 
\end{enumerate}
\end{proposition}

\begin{proof}
The statements can be deduced from Proposition \ref{prop:allstabilizergroups}, Lemmas \ref{lemma:zmod2stabilizer} and \ref{lemma:generallemma-groups}, and \cite[Proposition 5.19]{degaayfortman-nonarithmetic}. To show how this works, let us introduce some notation. Let $\tilde F \in \ca F_s(\RR)$ with image $F \in X_s(\RR)$. Let $f \in Y$ be the image of $[\tilde F] \in G(\RR) \setminus \ca F_s(\RR)$ under the isomorphism $G(\RR) \setminus \ca F_s(\RR) \xrightarrow{\sim} Y$ of Theorem \ref{th:realstableperiod}. Let $P\Gamma_f \subset P\Gamma$ be the stabilizer of $f \in Y$ in the group $P\Gamma$. Let $k \in \set{0,1,2}$ be the number of nodes of $F$, and write $k = 2a + b$, where $a$ is the number of pairs of complex conjugate nodes and $b$ the number of real nodes. Then the image $x \in \CC H^2$ of $\tilde F$ under \eqref{eq:stableperiodmapframed} lies on exactly $k$ distinct mutually orthogonal hyperplanes $H \in \ca H$, with $\ca H$ the set defined in \eqref{align:caH}. 
Since $F \in X_s(\RR)$, we have that $\tilde F \in \ca F_s^\alpha$ for some $\alpha \in P\mr A$. We get $x \in \RRH^n_\alpha$. 

Let $\ca H(x) \subset \ca H$ be the set of hyperplanes $H \in \ca H$ such that $x \in H$. 
Then $a$ equals the number of pairs of hyperplanes $H_1, H_2 \in \ca H(x)$ with $\alpha(H_1) = H_2$, and $b$ equals the number of hyperplanes $H \in \ca H(x)$ with $\alpha(H) = H$. 
Define $B_f \subset P\Gamma_f$ as the group generated by reflections in all $H \in \ca H(x)$ such that $\alpha(H) = H$. Consider the quotient map $p \colon \sqcup_{\alpha \in P\mr A} \RR H^2_\alpha \to Y$, let $\alpha_1, \dotsc, \alpha_\ell \in P\mr A$ be the elements such that $(x,\alpha_i) \sim (x,\alpha)$, and define $Y_f = \cup_{i = 1}^\ell p(\RR H^2_{\alpha_i}) \subset Y.$ The subgroup $P\Gamma_f \subset P\Gamma$ preserves the subset $Y_f \subset Y$ by \cite[Lemma 5.9]{degaayfortman-nonarithmetic}. Moreover, by \cite[Proposition 5.19.4]{degaayfortman-nonarithmetic}, there is an isometry between $B_f \sm Y_f$ and the union of $10^a$ copies of $\mathbb B^2(\RR)$. Let $\bb B$ be any one of these copies of $\BB^2(\RR)$, and define $$S_f \coloneqq \rm{Stab}_{P\Gamma_f/B_f}(\bb B),$$the stabilizer of $\bb B$ in the group $P\Gamma_f / B_f$. 
By construction of the orbifold $\overline{\mr M}_\RR$, see \cite[Propositions 5.1 \& 5.19]{degaayfortman-nonarithmetic}, the group $S_f$ is a representative of the isomorphism class of stabilizer groups of the orbifold $\overline{\mr M}_\RR$ at the moduli point $[F] \in \ca M_s(\RR)$ induced by $F$. Clearly, the stabilizer $G(\RR)_F \subset G(\RR)$ of $F \in X_s(\RR)$ yields the isomorphism class of stabilizer groups of the orbifold $\ca M_s(\RR)$ at the moduli point $[F] \in \ca M_s(\RR)$. In particular, we need to compare the isomorphism classes of the groups $S_f$ and $G(\RR)_F$. To do so, we claim that there is a canonical isomorphism \begin{align}\label{align:claim-frist}G(\RR)_F \xrightarrow{\sim} P\Gamma_f / G(x).\end{align}
Indeed, the actions of the groups $P\Gamma$ and $G(\RR)$ on $\ca F_s(\RR)$ commute, so we can apply Lemma \ref{lemma:generallemma-groups}. Recall that $\tilde F \in \ca F_s(\RR)$ has images $F \in X_s(\RR)$ and $f \in Y$, 
and that the map $\ca F_s(\RR) \to Y$ factors through an isomorphism $G(\RR) \sm \ca F_s(\RR) \xrightarrow{\sim} Y$. Moreover, if $P\Gamma_{\tilde F} \subset P\Gamma$ denotes the stabilizer of $\tilde F$ in $P\Gamma$, then 
$P\Gamma_{\tilde F} = G(x)$ by Lemma \ref{lemma:stabilizergroups}. As the group $G(\CC)$ acts freely on $\ca F_s$ by Lemma \ref{lemma:freeaction}, the group $G(\RR)$ acts freely on $\ca F_s(\RR)$. Thus, \eqref{align:claim-frist} follows from Lemma \ref{lemma:generallemma-groups}.

If $F$ has no nodes ($k = 0$), then $G(x)$ is trivial by \cite[Proposition 5.19.1]{degaayfortman-nonarithmetic}, and $P\Gamma_f = S_f$. Thus, $G(\RR)_F \cong S_f$ in view of \eqref{align:claim-frist}. 

If $F$ has only real nodes, then $B_f = G(x)$ and $P\Gamma_f/G(x) =S_f$. Thus, $G(\RR)_F \cong S_f$ by \eqref{align:claim-frist}. 

Finally, suppose that $F$ has a pair of complex nodes ($a = 1$ and $b = 0$). 
The zero set of $F$ defines a five-tuple $\underline{z}= (z_1, \dotsc, z_5) \in \PP^1(\CC)^5$, well-defined up to the $\PGL_2(\RR) \times \mf S_5$ action on $\PP^1(\CC)$, where $z_1 \in \PP^1(\RR)$ and $z_3 = \bar z_2  = z_5 = \bar z_4 \in \PP^1(\CC) \setminus \PP^1(\RR)$. Write $\underline{z} = (w, z, \bar z, z, \bar z)$ with $w \in \PP^1(\RR)$ and $z \in \PP^1(\CC) \setminus \PP^1(\RR)$. There is a unique $T \in \PGL_2(\RR)$ such that $T(w) = \infty$ and $T(z) = i$. This gives $T(\underline z) = (\infty, i , -i, i, -i)$. In particular, $F$ is unique up to isomorphism.

We have $G(x) \cong (\ZZ/10)^2$, and as there are no real nodes, $B_f$ is trivial. 
With respect to the isometry $\CC H^2(\CC) \xrightarrow{\sim} \BB^2(\CC)$ of \cite[Lemma 5.16]{degaayfortman-nonarithmetic}, the anti-holomorphic involutions $\alpha_j\colon\BB^2(\CC) \to \BB^2(\CC)$ induced by the $\alpha_j \in P\mr A$ are $(t_1, t_2) \mapsto (\bar t_2 \zeta^j , \bar t_1 \zeta^j)$, for $j \in \ZZ/10$. The fixed points sets are given by $\BB^2(\RR)_{\alpha_j}  = \{t_2 = \bar t_1 \zeta^j\}\subset \BB^2(\CC)$. The subgroup $E_f \subset G(x)$ that stabilizes $\BB^2(\RR)_{\alpha_j} $ is the group $E_f \cong \ZZ/10$ generated by the transformations $(t_1, t_2) \mapsto (\zeta t_1, \zeta^{-1} t_2)$. 
There is only one non-trivial $T \in \PGL_2(\RR)$ fixing $\infty$ and preserving $\{i, -i\} \subset \PP^1(\CC)$, and $T$ has order two, so $G(\RR)_F = \ZZ/2$. Therefore, by \eqref{align:claim-frist}, we have an exact sequence of groups $
0 \to G(x) \cong (\ZZ/10)^2 \to P\Gamma_f \to \ZZ/2 \to 0$ inducing an exact sequence $0 \to E_f \cong \ZZ/10 \to S_f \to \ZZ/2 \to 0$. 
By Proposition \ref{prop:allstabilizergroups} and Lemma \ref{lemma:zmod2stabilizer}, the proposition follows. 
\end{proof}
\subsection{The real moduli space as a hyperbolic triangle} \label{sec:triangle2} The goal of Section \ref{sec:triangle2} is to  prove Theorem \ref{theorem:lattice-thm}. 

The results in the above Sections \ref{sec-sec:1} and \ref{sec-sec:3} give the orbifold singularities of $\overline{\mr M}_{\RR}$ together with their stabilizer groups. In order to determine the hyperbolic orbifold structure of $\overline{\mr M}_{\RR}$, we also need to know the underlying topological space $\va{\ca M_s(\RR)}$ of $\overline{\mr M}_{\RR}$. The first observation is that $\ca M_s(\RR)$ is compact. Indeed, it is classical that the topological space $\ca M_s(\CC) = G(\CC) \setminus X_s(\CC)$, parametrizing complex stable binary quintics, is compact. This follows from the proper surjective map $\overline{M}_{0,5}(\CC)/ \mf S_5 \to \ca M_s(\CC)$ and the properness of the stack of stable five-pointed curves $\overline{M}_{0,5}$ \cite{Knudsen1983}, or from the fact that $\ca M_s(\CC)$ is homeomorphic to a compact ball quotient \cite{shimuratranscendental}. Moreover, the map $\ca M_s(\RR) \to \ca M_s(\CC)$ is proper, which proves the compactness of $\ca M_s(\RR)$. The second observation is that $\ca M_s(\RR)$ is connected, since $X_s(\RR)$ is obtained from the euclidean space $X(\RR) =  \{F \in \RR[x,y]\colon F \text{ homogeneous and} \deg(F) =  5\}$ by removing a subspace of codimension two. In the following lemma we generalize both of these observations. 
\begin{lemma} \label{lem:simplyconnected}
The moduli space $\ca M_s(\RR)$ is homeomorphic to a closed disc $\overline D \subset \RR^2$. 
\end{lemma}

\begin{proof}
The idea is to show that the following holds:
\begin{enumerate}
\item 
For each $j \in \{0,1,2\}$, the embedding $\mr M_j \hookrightarrow \overline{\mr M}_j \subset \ca M_s({\RR})$ of the connected component $\mr M_j$ of $\ca M_0(\RR)$ into its closure in $\ca M_s({\RR})$ is homeomorphic to the embedding $D \hookrightarrow \overline{D}$ of the open unit disc into the closed unit disc in $\RR^2$. 
\item We have $\ca M_s(\RR) = \overline{\mr M}_0  \cup \overline{\mr M}_1  \cup \overline{\mr M}_2$, and this gluing corresponds up to homeomorphism to the gluing of three closed discs $\overline{D}_i \subset \RR^2$ as in Figure~\ref{fig:triangle}. 
\end{enumerate}
To prove this, one considers the moduli spaces of real smooth (resp. stable) genus zero curves with five real marked points, as well as twists of this space. 
Define two anti-holomorphic involutions $\sigma_i\colon\PP^1(\CC)^5 \to \PP^1(\CC)^5$ by
$
\sigma_1(x_1, x_2, x_3, x_4, x_5) = (\bar x_1, \bar x_2, \bar x_3, \bar x_5, \bar x_4), 
$
and 
$
\sigma(x_1, x_2, x_3, x_4, x_5) = (\bar x_1, \bar x_3, \bar x_2, \bar x_5, \bar x_4). 
$
Then define 
$$
P_0^1(\RR) = P_0^{\sigma_1}, \white P_s^1(\RR) = P_1(\CC)^{\sigma_1}, \white
P_0^2(\RR) = P_0^{\sigma_2}, \white P_s^2(\RR) = P_1(\CC)^{\sigma_2}.
$$
It is clear that $ \mr M_0 = \PGL_2(\RR) \setminus P_0(\RR) / \mf S_5$. Similarly, we have:
\begin{align*}
   \mr M_1 = \PGL_2(\RR) \setminus P_0^1(\RR) / \mf S_3 \times \mf S_2 \quad \tn{ and } \quad  \mr M_2 = \PGL_2(\RR) \setminus P_0^2(\RR) / \mf S_2 \times \mf S_2.
    \end{align*}
  Moreover, we have $\overline{\mr M}_0 = \PGL_2(\RR) \setminus P_s(\RR) / \mf S_5$. We define 
\begin{align*}
\overline{\mr M}_1 = \PGL_2(\RR) \setminus P_s^1(\RR) / \mf S_3 \times \mf S_2, \quad \tn{ and } \quad 
\overline{\mr M}_2 = \PGL_2(\RR) \setminus P_s^2(\RR) / \mf S_2 \times \mf S_2.  
\end{align*}
Each $\overline{\mr M}_j$ is homeomorphic to a closed disc in $\RR^2$. Moreover, the natural maps $\overline{\mr M}_j \to \ca M_s(\RR) $ are closed embeddings of topological spaces, and one can check that the images glue to form $\ca M_s(\RR)$ in the prescribed way. We leave the details to the reader. 
\end{proof}

\begin{proof}[Proof of Theorem \ref{theorem:lattice-thm}]
To any closed two-dimensional orbifold $O$ one can associate a set of natural numbers $S_O = \{n_1, \dotsc, n_k; m_1, \dotsc, m_l\}$ by letting $k$ be the number of cone points of $X_O$, $l$ the number of corner reflectors, $n_i$ the order of the $i$-th cone point and $m_j$ the order of the $j$-th corner reflector (see \cite[Proposition 13.3.1]{Thurston80}). A closed two-dimensional orbifold $O$ is determined, up to orbifold-structure preserving homeomorphism, by its underlying space $X_O$ and the set $S_O$ \cite{Thurston80}. By Lemma \ref{lem:simplyconnected}, $\overline{\mr M}_\RR$ is homeomorphic to a closed disc in $\RR^2$. By Proposition \ref{prop:conesreflectors}, $\overline{\mr M}_\RR$ has no cone points and three corner reflectors whose angles are $\pi/3, \pi/5$ and $\pi/10$. This implies that $\overline{\mr M}_\RR$ and $\Delta_{3,5,10}$ are isomorphic as topological orbifolds. Consequently, the orbifold fundamental group of $\overline{\mr M}_\RR$ is abstractly isomorphic to the group $\Gamma_{3,5,10}$ defined in equation \eqref{PGAMMAR}.  

Let $\phi\colon\Gamma_{3,5,10} \hookrightarrow \text{PSL}_2(\RR)$ be \emph{any} embedding such that $X\coloneqq\phi\left(\Gamma_{3,5,10}\right) \setminus \RR H^2$ is a hyperbolic orbifold; we claim that there is a fundamental domain $\Delta$ for $X$ isometric to $\Delta_{3,5,10}$. To see this, consider the generator $a \in \Gamma_{3,5,10}$. Since $\phi(a)^2 = 1$, there exists a geodesic $L_1 \subset \RR H^2$ such that $\phi(a) \in \text{PSL}_2(\RR) = \text{Isom}(\RR H^2)$ is the reflection across $L_1$. Next, consider the generator $b \in \Gamma_{3,5,10}$. There exists a geodesic $L_2 \subset \RR H^2$ such that $\phi(b)$ is the reflection across $L_2$, and we have $L_2 \cap L_1 \neq \emptyset$. 
Let $x \in L_1 \cap L_2$. Then $\phi(a)\phi(b)$ is an element of order three that fixes $x$, hence is a rotation around $x$. Therefore, one of the angles between $L_1$ and $L_2$ must be $\pi/3$. Finally, we know that $\phi(c)$ is an element of order two in $\PSL_2(\RR)$, hence a reflection across a line $L_3$. As $L_3 \cap L_2 \neq \emptyset$ and $L_3 \cap L_1 \neq \emptyset$, 
the three geodesics $L_i \subset \RR H^2$ enclose a hyperbolic triangle. As the orders of the three elements $\phi(a)\phi(b)$, $\phi(a)\phi(c)$ and $\phi(b)\phi(c)$ are respectively three, five and ten, 
the three interior angles of the triangle are $\pi/3$, $\pi/5$ and $\pi/10$. Thus, $X$ is isometric to $\Delta_{3,5,10}$. 

Consequently, $P\Gamma_\RR \sm \RRH^2$ is isometric to $\Gamma_{3,5,10} \sm \RRH^2$. It follows that the lattices $P\Gamma_\RR$ and $\Gamma_{3,5,10}$ are conjugate in $\rm{PO}(2,1)$, see e.g.\ \cite[Lemma 1]{ratcliffe}. 
\end{proof}

\section{The monodromy groups} \label{sec:monodromy} In this section, we describe the monodromy group $P\Gamma$ attached to the moduli space $X_0(\CC)$, as well as the groups $P\Gamma_\alpha$ appearing in Proposition \ref{prop:realsmoothperiods}. As for the lattice $(\Lambda, \mf h)$ (see (\ref{eq:hermitianformonbinaryquinticlattice})), we have:
\begin{theorem}[Shimura] \label{th:calculatemonodromyshimura}
There is an isomorphism of hermitian $\OO_K$-lattices
 $$\left( \Lambda, \mf h \right) \cong \left( \OO_K^3, \textnormal{diag}\left(- \lambda, 1,1\right) \right), \quad \lambda = \zeta_5 + \zeta_5^{-1} = \frac{\sqrt 5 - 1}{2}.$$
\end{theorem}
\begin{proof}
See \cite[Section 6]{shimuratranscendental} as well as item (5) in the table on page 1 of \emph{loc.\ cit.}
\end{proof}
Write $\Lambda = \OO_K^3$ and $\mf h = \textnormal{diag}(- \lambda, 1, 1)$. Consider the $\bb F_5$-vector space $W = \Lambda / (1 - \zeta_5)\Lambda$ equipped with the quadratic form $q = \mf h \bmod \theta$. 
Define three anti-isometric involutions as follows:
\begin{equation} \label{chi}\begin{split} 
    \alpha_0\colon & (x_0,x_1,x_2) \mapsto (\bar x_0, \;\;\;\bar x_1,\;\;\bar x_2) \\
    \alpha_1 \colon & (x_0,x_1,x_2) \mapsto (\bar x_0,  - \bar x_1, \;\;\bar x_2) \\
    \alpha_2 \colon & (x_0,x_1,x_2) \mapsto (\bar x_0, -\bar x_1, -\bar x_2).
    \end{split}
\end{equation}
\begin{lemma} \label{exactlyconjugate}
An anti-unitary involution of $\Lambda$ is $\Gamma$-conjugate to exactly one of the $\pm \alpha_j$. 
In particular, $\va{P\Gamma \sm P\mr A} = 3$ and the $\alpha_i$ defined in (\ref{chi}) form a set of representatives for $P\Gamma \sm P\mr A$. 
\end{lemma}
\begin{proof}
For isometries $\alpha \colon W \to W$, the dimension and determinant of the fixed space $(W^\alpha, q|_{W^\alpha})$ are conjugacy-invariant. Using this, one can show that 
the elements $\pm \alpha_i$ are pairwise non $\Gamma$-conjugate. Moreover, $\va{P\Gamma \sm P\mr A} = \va{\pi_0(\ca M_0(\RR))} = 3$ by Proposition \ref{prop:realsmoothperiods} and Theorem \ref{th:calculatemonodromyshimura}.
\end{proof}
Define $$\theta = \zeta_5 - \zeta_5^{-1} \in \OO_K.$$ Note that $|\theta|^2 = \frac{\sqrt{5} + 5}{2}$. The fixed lattices of the anti-unitary involutions $\alpha_i$ defined in \eqref{chi} are
\begin{align}\label{eq:fixedlattices} 
\begin{split}
\Lambda^{\alpha_0} &= \ZZ[\lambda] \oplus \ZZ[\lambda] \oplus \ZZ[\lambda], \\
\Lambda^{\alpha_1} &= \ZZ[\lambda] \oplus \theta\ZZ[\lambda] \oplus \ZZ[\lambda],\\
\Lambda^{\alpha_2} &= \ZZ[\lambda] \oplus \theta\ZZ[\lambda] \oplus  \theta\ZZ[\lambda],
\end{split}
\end{align}
where \begin{align} \label{align:lambda-def} \lambda = \zeta_5 + \zeta_5^{-1} = \frac{\sqrt 5 - 1}{2}.\end{align}
Restricting $\mf h$ to the $\Lambda^{\alpha_j}$ yields quadratic forms $q_0$, $q_1$ and $q_2$ on $\ZZ[\lambda]^3$ defined as follows:
\begin{equation} \label{eq:explicitquadraticforms}
\begin{split}
 q_0(x_0, x_1, x_2) & = - \lambda x_0^2 + x_1^2 + x_2^2, 
\\
 q_1(x_0, x_1, x_2) & =  - \lambda x_0^2 + \left(\frac{\sqrt{5} + 5}{2}\right) \cdot x_1^2 + x_2^2, 
\\
 q_2(x_0, x_1, x_2) & = - \lambda x_0^2 + \left(\frac{\sqrt{5} + 5}{2}\right) \cdot x_1^2 + \left(\frac{\sqrt{5} + 5}{2}\right) \cdot x_2^2, 
 \end{split}
\end{equation}
We consider $\ZZ[\lambda]$ (with $\lambda$ as in \eqref{align:lambda-def}) as a subring of $\RR$ via the embedding that sends $\lambda$ to a positive element. 
\begin{theorem} \label{th:explicitquadratic}
Consider the quadratic forms $q_j$ defined in (\ref{eq:explicitquadraticforms}), with $\lambda$ as in \eqref{align:lambda-def}. There is a union of geodesic subspaces $\mr H_j \subset \RR H^2$ $\left( j \in \{0,1,2\} \right)$ and an isomorphism of hyperbolic orbifolds
\begin{equation} 
    \ca M_0(\RR) \cong \coprod_{j = 0}^2 \textnormal{PO}(q_j,\ZZ[\lambda]) \setminus \left(\RR H^2 - \mr H_j \right). 
\end{equation}
\end{theorem}

\begin{proof}
By Proposition \ref{prop:realsmoothperiods} and Lemma \ref{exactlyconjugate}, we obtain an isomorphism $$\ca M_0(\RR) \cong \coprod_{j = 0}^2 P\Gamma_{\alpha_j} \setminus (\RR H^2_{\alpha_j} - \mr H).$$ 
Remark that $P\Gamma_{\alpha_j}=N_{P\Gamma}(\alpha_j)$ for the normalizer $N_{P\Gamma}(\alpha_j)$ of $\alpha_j$ in $P\Gamma$. If $h_j$ denotes the restriction of $\mf h$ to $\Lambda^{\alpha_j}$, there is a natural embedding
$$
\iota_j \colon    N_{P\Gamma}(\alpha_j) \hookrightarrow \textnormal{PO}(\Lambda^{\alpha_j},h_j, \ZZ[\lambda]).$$ We claim that $\iota_j$ is an isomorphism. This holds, because the natural homomorphism
$
\pi_j \colon N_\Gamma(\alpha_j) \to \rm O(\Lambda^{\alpha_j}, h_j)
$
is surjective, where $N_\Gamma(\alpha_j) = \{g \in \Gamma\colon g \circ \alpha_j = \alpha_j \circ g \}$ is the normalizer of $\alpha_j$ in $\Gamma$. The surjectivity of $\pi_j$ follows in turn from the equality
\begin{align} \label{equality:equality}
    \Lambda =  \ca O_K \cdot  \Lambda^{\alpha_j} + \ca O_K \cdot \theta \left(\Lambda^{\alpha_j}\right)^\vee \subset K^3,
\end{align}
and \eqref{equality:equality} follows from \eqref{eq:fixedlattices}. Since $\textnormal{PO}(\Lambda^{\alpha_j},h_j, \ZZ[\lambda]) = \textnormal{PO}(q_j,\ZZ[\lambda])$, we are done. 
\end{proof}

\begin{proof}[Proof of Theorem \ref{th:theorem02}]
This follows from Corollary \ref{cor:theorem2}, Theorem \ref{theorem:lattice-thm} and Theorem \ref{th:explicitquadratic}.  
\end{proof}

\begin{proof}[Proof of Theorem \ref{theorem:angles}]
In \cite{aperyyoshida-pentagonalstructure}, Apéry and Yoshida proved that $\overline{\mr M}_0$ is the hyperbolic triangle of angles $\pi/2, \pi/4$ and $\pi/5$. 
As the two hyperplanes in Figure \ref{fig:triangle} intersect orthogonally, this implies that the bottom angle of the triangle $\overline{\mr M}_0$ in Figure \ref{fig:triangle} (i.e.~its angle at $(0,-1,\infty,\infty,1)$) equals $\pi/2$, and that the angle of $\overline{\mr M}_0$ at $(0,-1,-1,\infty,\infty)$ equals $\pi/4$. One deduces that the left angle of $\overline{\mr M}_1$ is $\pi/2$, and that the angle of $\overline{\mr M}_2$ at $(0,-1,-1,\infty,\infty)$ equals $\pi/4$. 

For a hyperbolic triangle with angles $\alpha, \beta, \gamma$ and sides $a,b,c$ such that $a$ is the side opposite to $\alpha$, $b$ the side opposite to $\beta$ and $\gamma$ the side opposite to $c$, one has the hyperbolic law of cosines 
\begin{align}\label{align:cosines}
\cosh(c) = \frac{\cos(\alpha)\cos(\beta)+\cos(\gamma)}{\sin(\alpha)\sin(\beta)}.
\end{align}
Applying \eqref{align:cosines} to the triangles $\overline{\mr M}_0$ and $\overline{\mr M}_\RR$, one can calculate the length of the side of $\overline{\mr M}_2$ that connects $(0,-1,-1,\infty,\infty)$ and $(\infty, i, i, -i, -i)$. Applying \eqref{align:cosines} again, it follows that the angle of $\overline{\mr M}_2$ at the point $(0,-i,\infty,\infty,i)$ is $\pi/2$. Thus, the angle of $\overline{\mr M}_1$ at $(0,-i,\infty,\infty,i)$ is also $\pi/2$.
\end{proof}
%
%

\section{Non-arithmetic lattices in the projective orthogonal group} \label{sec:announced}

In a previous paper we proved a result, see \cite[Theorem 1.8]{degaayfortman-nonarithmetic}, that has the following consequence. 
For $n \geq 2$, define 
\[
\mr L_{\zeta_5}^n(\lambda) \coloneqq \left(\QQ(\zeta_5), \ZZ[\zeta_5]^{n,1}_\lambda \right), \quad \quad \lambda = \zeta_5 + \zeta_5^{-1} = (\sqrt 5 - 1)/2. 
\]
Here, $ \ZZ[\zeta_5]^{n,1}_\lambda$ is the free $\ZZ[\zeta_5]$-module of rank $n+1$ equipped with the hermitian
form $h$ defined as $h(x,y) = -\lambda \cdot x_0 \bar y_0 + \cdots + x_n \bar y_n$. Then $\mr L_{\zeta_5}^n(\lambda) $ is a hermitian lattice of rank $n+1$ in the sense of \cite[Definition 2.2]{degaayfortman-nonarithmetic} (indeed, this follows from \cite[Example 2.12]{degaayfortman-nonarithmetic}). 
For each $n \geq 2$, perform the gluing construction of \cite[Definition 1.1]{degaayfortman-nonarithmetic} to associate to the hermitian lattice $\mr L_{\zeta_5}^n(\lambda)$ a topological space 
$
M( \mr L_{\zeta_5}^n(\lambda)).
$
By \cite[Theorem 1.2]{degaayfortman-nonarithmetic}, there exists a canonical real hyperbolic orbifold structure on $M( \mr L_{\zeta_5}^n(\lambda)  )$ such that each connected component of $M( \mr L_{\zeta_5}^n(\lambda)  )$ is isomorphic to the quotient of real hyperbolic $n$-space $\RRH^n$ by a lattice in $\rm{PO}(n,1)$.  Define an anti-unitary involution $\alpha_0 \colon \ZZ[\zeta_5]^{n,1}_\lambda  \to \ZZ[\zeta_5]^{n,1}_\lambda $ by $\alpha_0(x) = \bar x$, let $$M\left( \mr L_{\zeta_5}^n(\lambda)  , \alpha_0 \right) \subset M\left( \mr L_{\zeta_5}^n(\lambda)  \right)$$ be the connected component that contains the image of the natural map $\RRH^n_{\alpha_0} \to M( \mr L_{\zeta_5}^n(\lambda)  )$, and let 
\[
\Gamma_{\zeta_5}^n(\lambda) \subset \rm{PO}(n,1)
\]
be a lattice such that $M( \mr L_{\zeta_5}^n(\lambda)  , \alpha_0)  \cong \Gamma_{\zeta_5}^n(\lambda)  \sm \RRH^n$ (compare \cite[page 7]{degaayfortman-nonarithmetic}). 

By combining \cite[Theorem 1.8]{degaayfortman-nonarithmetic} with the main results of this paper, one can prove the following result. 

\begin{theorem}
For each $n \geq 2$, the lattice $\Gamma_{\zeta_5}^n(\lambda) \subset \rm{PO}(n,1)$ is non-arithmetic. 
\end{theorem}
\begin{proof}
Write $\ZZ[\zeta_5]^{2,1}_\lambda = \Lambda$; this abuse of notation is harmless in view of Theorem \ref{th:calculatemonodromyshimura}. Let $\Gamma = \Aut(\Lambda)$. 
By Lemma \ref{exactlyconjugate}, 
an anti-unitary involution of $\Lambda$ is $\Gamma$-conjugate to exactly one of the involutions $\pm \alpha_j$ defined in equation \eqref{chi}. We can therefore apply \cite[Theorem 1.8]{degaayfortman-nonarithmetic}, which implies that $\Gamma_{\zeta_5}^n(\lambda) \subset \rm{PO}(n,1)$ is non-arithmetic for each $n \geq 2$ provided that $\Gamma_{\zeta_5}^2(\lambda) \subset \rm{PO}(2,1)$ is non-arithmetic. In other words, we are reduced to the case $n=2$. By Theorem \ref{th:calculatemonodromyshimura}, the lattice $\Gamma_{\zeta_5}^2(\lambda) \subset \rm{PO}(2,1)$ is conjugate to the lattice $P\Gamma_\RR \subset \rm{PO}(2,1)$ defined in Corollary \ref{cor:theorem2}. Moreover, by Theorem \ref{theorem:lattice-thm}, the lattice $P\Gamma_\RR$ is conjugate to the lattice $\Gamma_{3,5,10} \subset \rm{PO}(2,1)$ defined in equation \eqref{PGAMMAR}. Finally, by Takeuchi's classification of arithmetic triangle groups, see \cite{takeuchi}, the subgroup $\Gamma_{3,5,10} \subset \rm{PO}(2,1)$ is non-arithmetic. Thus, $\Gamma_{\zeta_5}^2(\lambda)$ is non-arithmetic, and the theorem follows. 
\end{proof}

\begin{acknowledgements}
This project was carried out partly at the ENS in Paris and partly at the Leibniz University in Hannover. I thank my former PhD advisor Olivier Benoist for his guidance and support. I thank Romain Branchereau, Samuel Bronstein, Nicolas Tholozan and Frans Oort for useful discussions. I thank Nicolas Bergeron for pointing me to Takeuchi's paper on hyperbolic triangle groups. 

I would like to thank the referee for his or her valuable comments on this paper. 

This project has received funding from the European Union's Horizon 2020 research and innovation programme under the Marie Sk\l{}odowska-Curie grant agreement N\textsuperscript{\underline{o}} 754362 \img{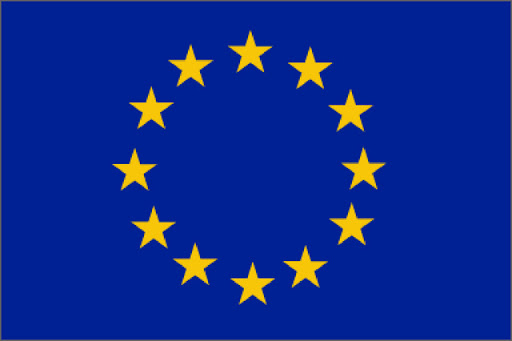} 
and from the European Research Council (ERC) under the European Union’s Horizon 2020 research and innovation programme under grant agreement N\textsuperscript{\underline{o}} 948066 (ERC-StG RationAlgic). 
\end{acknowledgements}

\printbibliography

\end{document}

%% file: macros.tex
\newcommand*{\img}[1]{%
    \raisebox{-.02\baselineskip}{%
        \includegraphics[
        height=\baselineskip,
        width=\baselineskip,
        keepaspectratio,
        ]{#1}%
    }%
}

\newcommand{\PSL}{\textnormal{PSL}}

\newcommand{\set}[1]{\left\{ #1 \right\}}
\newcommand{\va}[1]{\left| #1 \right|}

\newcommand{\Id}{\textnormal{Id}}

\newcommand{\white}{\;\;\;}

\newcommand{\Stab}{\textnormal{Stab}}
\newcommand{\sm}{\setminus}

\newcommand{\CC}{\mathbb{C}}

\newcommand{\BB}{\mathbb{B}}

\newcommand{\OO}{\mathcal{O}}

\newcommand{\VV}{\mathbb{V}}

\newcommand{\QQ}{\mathbb{Q}}
\newcommand{\PP}{\mathbb{P}}

\newcommand{\PGL}{\textnormal{PGL}}

\newcommand{\GL}{\textnormal{GL}}

\newcommand{\SL}{\textnormal{SL}}

\newcommand{\RR}{\mathbb{R}}

\newcommand{\ZZ}{\mathbb{Z}}

\newcommand{\CCH}{\bb C H}
\newcommand{\RRH}{\bb R H}

\newcommand{\Ker}{\textnormal{Ker}}

\newcommand{\id}{\textnormal{id}}

\newcommand{\Hom}{\textnormal{Hom}}

\newcommand{\Gal}{\textnormal{Gal}}

\newcommand{\End}{\textnormal{End}}

\newcommand{\Aut}{\textnormal{Aut}}

\newcommand{\ca}[1]{{\mathcal{#1}}}
\newcommand{\bb}[1]{{\mathbb{#1}}}

\newcommand{\mr}[1]{{\mathscr{#1}}}
\newcommand{\mf}[1]{{\mathfrak{#1}}}
\newcommand{\tn}[1]{{\textnormal{#1}}}
\let\rm\relax 
\newcommand{\rm}[1]{{\mathrm{#1}}}

\DeclareSymbolFont{bbm}{U}{bbm}{m}{n}
\DeclareSymbolFontAlphabet{\mathbbm}{bbm}
